\documentclass[12pt]{amsart}
\usepackage{mathrsfs}
\usepackage{mathdg2015-SPP2025}

\usepackage[utf8]{inputenc}
\usepackage{amsmath}
\usepackage{amsfonts}
\usepackage{amssymb}
\usepackage{amsthm}
\usepackage{float}
\usepackage{mathtools}
\usepackage{csquotes}
\usepackage{pdfpages}
\usepackage{enumitem}
\usepackage[margin=3cm]{geometry}

\usepackage{tikz-cd}

\usepackage{mathtools}
\usepackage{leftindex}
\usepackage{slashed}
\usepackage{circuitikz}
\usepackage{hyperref}

\newcommand{\DN}{\mathcal{N}}

\newcommand{\bbf}{\mathrm{bf}}
\newcommand{\phif}{\phi\mathrm{f}}

\newcommand{\s}{\mathrm{s}}

\newcommand{\dsX}{{\partial_{\s} X}}

\newcommand{\textBC}{{\upshape\textsc{bc}}}
\newcommand{\BC}{\ensuremath{\text{\resizebox{!}{4pt}{\tiny\upshape BC}}}}

\newcommand{\dBCX}{{\partial_\BC X}}

\newcommand{\dsBCX}{{\partial_{\s,\BC} X}}

\newcommand{\Dnu}{D_\nu}
\newcommand{\phisigma}{{}^\phi\sigma} 
\newcommand{\phiT}{{}^\phi T} 

\newcommand{\cl}{{\mathrm{cl}}}
\newcommand{\Sph}{\mathbb{S}}
\newcommand{\Cdotinf}{\dot{C}^\infty}
\newcommand{\bdX}{\partial X}
\newcommand{\Uwidetilde}{\widetilde{U}}

\renewcommand{\dX}{{\partial X}}

\newcommand{\codim}{\operatorname{codim}}

\newtheorem{assumption}[theorem]{Assumption}

\newcommand{\dirac}{\slashed{\partial}}
\newcommand\restr[2]{\ensuremath{\left.#1\right|_{#2}}}

\newcommand{\cU}{\calU}

\makeatletter
\@namedef{subjclassname@2020}{%
  \textup{2020} Mathematics Subject Classification}
\makeatother

\graphicspath{{figures/}}

\title{Analysis on fibred cusp spaces}
\author{Daniel Grieser}
\address{Institut f\"ur Mathematik, Carl von Ossietzky Universit\"at Oldenburg}
\email{daniel.grieser@uol.de}
\author{\'Alvaro S\'anchez-Hern\'andez}
\address{Institut f\"ur Analysis, Leibniz Universit\"at Hannover}
\email{ash@math.uni-hannover.de}
\author{Boris Vertman}
\address{Institut f\"ur Mathematik, Carl von Ossietzky Universit\"at Oldenburg}
\email{boris.vertman@uol.de}
\thanks{DG and BV were partially supported by the DFG Priority Programme 2026 \lq Geometry at Infinity\rq.}

\date{\today}


\begin{document}

\pagestyle{plain}

\begin{abstract}
We give a survey of analytic and geometric results on `fibred cusp spaces', a large class of non-compact Riemannian manifolds
which include the regular parts of singular spaces with incomplete cusp singularities as well as complete spaces with asymptotically hyperbolic cusp or asymptotically Euclidean structures at infinity. These results cover topics in spectral geometry, in particular analytic torsion and index theory, and boundary value problems. The underlying tools include a careful microlocal analysis of the resolvent and the heat kernel. We include an exposition of the geometric and analytic foundations and sketch the ideas of the proofs of the main theorems. Special emphasis is put on the common features of and the differences between the incomplete and various kinds of complete settings.    
\end{abstract}

\maketitle

\phantomsection

\keywords

\tableofcontents

\section{Introduction}

\bigskip

The word cusp(idal) is used in various areas of mathematics for sharply pointed objects. In geometry, two common types of cusps are the incomplete cusp and the hyperbolic cusp, shown in Figure \ref{fig:examples} (left and center). An important feature is that the  metric near the singularity of the incomplete cusp is conformal to the metric near infinity of the hyperbolic cusp. It turns out that this conformal class has a natural  generalization which includes the metric on the large end of a cone, shown on the right of Figure \ref{fig:examples}. In addition, metrics in this class have asymptotic smoothness properties near the respective singularity/infinity which can be expressed in terms of a compactification of the underlying space to a manifold with boundary, where the boundary corresponds to the singularity or to infinity. We call such metrics  $c-\phi$-metrics, where $c$ refers to the exponent appearing in the conformal factor and $\phi$ to a fibration of the boundary inherent to the geometry. Also, we call a manifold with fibred boundary, equipped with such a metric on its interior, a fibred cusp space. 
Such spaces arise in many widely different contexts. 

Over the last three decades there has been a lot of progress on various questions related to the analysis on such spaces, including Fredholmness conditions on geometric elliptic operators, heat kernel asymptotics, index theory and more. The approach taken in many of these works, and in this paper, is via microlocal analysis, and more specifically geometric singular microlocal analysis, as originally introduced by R. Melrose. The class of $\phi$-manifolds, which is the basis for our setting, was introduced by Mazzeo and Melrose in \cite{MM98} (they are called manifolds with fibred boundary there).

In this article we give an overview of some of these works. The selection is certainly incomplete and influenced by the research interests of the authors. 

\medskip

Outline: In Section \ref{sec:fib cusp} we define $\phi$-manifolds, $c-\phi$-metrics and fibred cusp spaces, give examples and mention some generalizations of them.  The basic analytic framework, the so-called $\phi$-calculus, is introduced in Section \ref{sec:phi-calculus}. We then discuss the precise construction of the resolvent and of the heat kernel of Laplace type operators. These are used to analyse analytic torsion and index theory in Sections \ref{sec:torsion} and \ref{sec:index}. In Section \ref{sec:calderon} we generalize the setting somewhat, by allowing the fibres of a fibred cusp space to have boundary, and discuss the Calderón projector (for general elliptic operators) and the Dirichlet-Neumann operator (for the Laplacian) for boundary value problems with respect to that boundary. We conclude in Section \ref{sec:more} with a few remarks on other research directions and on open problems.
Throughout we use the language of manifolds with corners and blow-ups, which we explain in the Appendix.

\begin{figure}[!ht]
\centering
\resizebox{1\textwidth}{!}{%
\begin{circuitikz}
\tikzstyle{every node}=[font=\LARGE]
\draw [short] (5,5) .. controls (4.5,9) and (3.25,8.75) .. (1.25,10);
\draw [short] (5,5) .. controls (5.75,9) and (7,8.75) .. (8.75,10);
\draw [short] (4.6,7) .. controls (4.75,6.75) and (5.25,6.75) .. (5.55,7);
\draw [dashed] (5.55,7) .. controls (5.25,7.25) and (4.75,7.25) .. (4.6,7);
\draw [short] (2.65,9.25) .. controls (3,8.75) and (7,8.75) .. (7.5,9.25);
\draw [dashed] (7.5,9.25) .. controls (6.75,9.75) and (3.25,9.75) .. (2.65,9.25);
\draw [short] (1.25,10) .. controls (-3.75,14) and (0.5,18.5) .. (5,18.25);
\draw [short] (8.75,10) .. controls (14,14) and (9.75,18.5) .. (5,18.25);
\draw [short] (-0.95,13.25) .. controls (1.25,12.25) and (8.75,12.25) .. (11.1,13.25);
\draw [dashed] (11.1,13.25) .. controls (8.75,14) and (1.25,14) .. (-0.95,13.25);
\node [font=\Huge] at (5,1.5) {\fontsize{40}{28}\selectfont $dx^2 + x^4dz^2$};
\draw [line width=1pt, ->, >=Stealth] (2.5,12) .. controls (3.25,11.75) and (6.75,11.75) .. (7.5,12) ;
\node [font=\Huge] at (5,11.25) {\fontsize{30}{28}\selectfont $z$};
\draw [line width=1pt, ->, >=Stealth] (11.5,5.25) -- (11.5,8.75);
\node [font=\Huge] at (12,7) {\fontsize{30}{28}\selectfont $x$};
\draw [short] (20,5) -- (20,16.25);
\draw [short] (26.25,5) -- (26.25,16.25);
\draw [short] (20,5) .. controls (21.25,7.5) and (25,7.5) .. (26.25,5);
\draw [ color={rgb,255:red,247; green,8; blue,8}, line width=2pt, dashed] (20,18.75) -- (26.25,18.75);
\draw [dashed] (20,16.25) -- (20,18.5);
\draw [dashed] (26.25,16.25) -- (26.25,18.75);
\draw [line width=1pt, ->, >=Stealth] (18.75,18.25) -- (18.75,15);
\draw [line width=1pt, ->, >=Stealth] (21,13.75) -- (25.5,13.75);
\node [font=\Huge] at (22.5,1.5) {\fontsize{40}{28}\selectfont $\frac{dx^2}{x^2}+x^2dz^2$};
\node [font=\Huge] at (23.25,13.25) {\fontsize{30}{28}\selectfont $z$};
\node [font=\Huge] at (18.25,16.75) {\fontsize{30}{28}\selectfont $x$};
\draw [short] (35,16.25) .. controls (37.5,9.5) and (40,5) .. (41.25,5);
\draw [short] (41.25,5) .. controls (42.5,5) and (45,9.5) .. (47.5,16.25);
\draw [dashed] (35,16.25) -- (33.75,19.25);
\draw [dashed] (47.5,16.25) -- (48.75,19.25);
\draw  (41.25,16.25) ellipse (6.25cm and 1cm);
\draw [ color={rgb,255:red,247; green,8; blue,8} , line width=2pt , dashed] (41.25,19.25) ellipse (7.5cm and 0.5cm);
\draw [line width=1pt, ->, >=Stealth] (33,18.75) -- (33,14.25);
\draw [line width=1pt, ->, >=Stealth] (38,14.5) .. controls (39.25,14.25) and (43.25,14.25) .. (44.5,14.5) ;
\node [font=\Huge] at (32.55,16.5) {\fontsize{30}{28}\selectfont $x$};
\node [font=\Huge] at (41.25,13.75) {\fontsize{30}{28}\selectfont $y$};
\node [font=\Huge] at (41.25,1.5) {\fontsize{40}{28}\selectfont $\frac{dx^2}{x^4}+\frac{dy^2}{x^2}$};
\draw [line width=1pt, ->, >=Stealth] (18.75,8) -- (18.75,11.25);
\draw [line width=1pt, ->, >=Stealth] (33.25,8) -- (33.25,11.25);
\node [font=\Huge] at (18.1,9.5) {\fontsize{30}{28}\selectfont $w$};
\node [font=\Huge] at (32.75,9.5) {\fontsize{30}{28}\selectfont $r$};
\draw [ color={rgb,255:red,247; green,8; blue,8} , fill={rgb,255:red,247; green,8; blue,8}] (5,5) circle (0.25cm);
\end{circuitikz}
}%

\caption{Three instances of fibred cusp spaces. On the left, we have an incomplete cusp with fibre $\mathbb{S}^1$ (the corresponding space $X$ has the singular point replaced by a circle, given in coordinates by $x=0$). In the middle, we start with a fundamental domain $\text{SL}(2,\mathbb{R})/\text{SL}(2,\mathbb{Z})$ in the upper half plane model of hyperbolic space with the hyperbolic metric $\frac{dz^2 + dw^2}{w^2}$ and compactify at large $w$ by introducing the coordinate $x=\frac{1}{w}$ and adding $x=0$ (also, one should cut off at finite $w$ and add a smooth cap). On the right, we compactify the large end of a cone $dr^2+r^2dy^2$ by setting $x=\frac{1}{r}$ and adding $x=0$ (and smooth out at $r=0$); an example is $\R^n$ with the Euclidean metric, compactified by adding a sphere at $\infty$. In this case $B=\mathbb{S}^{n-1}$ and $F$ is a point.}
\label{fig:examples}
\end{figure}

\section{Fibred cusp spaces} \label{sec:fib cusp}

\subsection{Definitions}\label{subsec:def fib cusp}

By definition, a \textbf{$\phi$-manifold} is a compact manifold $X$ with boundary,  equipped with a fibration
\begin{equation}
\label{eqn:def phi-manifold}
\bdX \stackrel{\phi}\to B
\end{equation}
for a closed manifold  $B$, the base. 
To simplify the exposition we assume that the boundary is connected; however, see Remarks \ref{rem:connected}, \ref{rem:connected2}. We denote the fibre by $F$ and let $n=\dim X$, $b=\dim B$, $f=\dim F$.
For simplicity we assume that a trivialization of a neighborhood of the boundary
$$\calU \cong [0,\eps) \times \bdX$$
is given and that $\phi$ is extended to $\calU$. We denote local coordinates on $[0,\eps)$ by $x$ (so $\bdX=\{x=0\}$), on $B$ (and then pulled back to $\calU$) by $y=(y_1,\dots,y_b)$ and on $F$ by $z=(z_1,\dots,z_f)$.
Our main object of study are \textbf{general $c-\phi$-metrics} for $c\in\Z$, which by definition are Riemannian metrics on the interior $\interior{X}$ which on $\interior{\calU}$ have the form
$$ g_c = x^{2c} g_\phi, $$
where $g_\phi$ is a \textbf{general $\phi$-metric}, which means that in local coordinates it is a quadratic form expressible via the 1-forms
\begin{equation}
\label{eqn:phi-forms}
\frac{dx}{x^2}\,,\quad \frac{dy_i}x\,,\quad dz_j,
\end{equation}
($i\in\{1,\dots,b\}$ and $j\in\{1,\dots,f\}$) 
with a coefficient matrix that is smooth and positive definite up to the boundary $\{x=0\}$.

Invariantly, this may be expressed by saying that $g_\phi$ is a metric on a rescaled tangent bundle, $\phiT X$, which is defined as follows.
Define the Lie algebra of vector fields $\calV_\phi(X)$ as consisting of those smooth vector fields on $X$ which are tangent to the boundary, are tangent to the fibres of $\phi$ there, and satisfy $Vx=O(x^2)$.\footnote{When replacing $x$ by a different boundary defining function $\xtilde$, i.e. $\xtilde=ax$ with $a>0$ smooth, the last condition is unchanged iff $a$ is constant along the fibres of $\bdX$. Therefore, the data needed to fix $\calV_\phi(X)$, and hence the notion of general $\phi$-metric, is $\phi$ in \eqref{eqn:def phi-manifold} plus a choice of $x$ modulo this equivalence.} 
In local coordinates near the boundary this space is spanned, over $\Cinf(X)$, by the $n$ vector fields
\begin{equation}
\label{eqn:phi-vector fields}
x^2\partial_x\,,\quad x\partial_{y_i}\,,\quad \partial_{z_j}.
\end{equation}
Being a locally free $\Cinf(X)$-module of rank $n$, the space $\calV_\phi(X)$ defines a rank $n$ vector bundle (by the Serre-Swan theorem)
\begin{equation}
\label{eqn:def phiTX}
\phiT X \to X,
\end{equation}
together with a vector bundle map $\iota:\phiT X\to TX$, which is an isomorphism over the interior, 
by the requirement that the inclusion $\calV_\phi(X)\to\Cinf(X,TX)$ lifts to an isomorphism $\calV_\phi(X)\to\Cinf(X,\phiT X)$ under the induced map 
$\iota:\Cinf(X,\phiT X)\to\Cinf(X,TX)$. In particular, \eqref{eqn:phi-vector fields} is a local basis for $\phiT X$, and \eqref{eqn:phi-forms} is a local basis for $\phiT^*X$.

A \textbf{fibred cusp space} is a $\phi$-manifold equipped with a general $c-\phi$-metric, for some $\phi$ and $c$. Usually we take $c\in\{0,1,2\}$, see the examples below.

For some of the theorems stated in this paper additional conditions on the general $\phi$-metric $g_\phi$ must be imposed:
We call $g_\phi$ a \textbf{$\phi$-metric} if on $\calU$ it has the form
\begin{equation}
\label{eqn:phi metric}
\frac{dx^2}{x^4} + \frac{\phi^*g_B}{x^2} + h,
\end{equation} 
where $g_B$ is a family of metrics on $B$ smoothly depending on $x\in[0,\eps)$ and $h$ is a smooth symmetric 2-tensor on $TX$ which restricts to be positive definite on $TF$.\footnote{In coordinates this means that the coefficients of mixed terms $dx\,dy_i$, $dx\,dz_j$, $dy_i\,dz_j$ are smooth rather than just $x^{-3}$, $x^{-2}$, $x^{-1}$ times smooth, respectively, that the coefficient of $\frac{dx^2}{x^4}$ is one modulo $x^4$ and that the coefficients of $\frac{dy_i\,dy_{i'}}{x^2}$ are independent of $z$ modulo $x^2$.} We call $g_\phi$ a \textbf{product type} $\phi$-metric if $g_B$ and $h$ are constant in $x$.
Sometimes we need to require in addition that $\phi:(\bdX,\phi^*g_B + h)\to (B,g_B)$ be a \textbf{Riemannian submersion}.

In the literature  various names are used for $c-\phi$-metrics, for example:
\begin{itemize}
\item $c=0$: fibred boundary metric, $\phi$-metric, scattering metric (if $F$ is a point)
\item $c=1$: (generalized)
fibred cusp metric, d-metric, fibred hyperbolic cusp metric
\item $c=2$: incomplete fibred cusp metric
\end{itemize}
If the base $B$ is a point then 
$c-\phi$-metrics are often called cusp metrics.

\begin{remark}[Connectedness I]
\label{rem:connected}
If $\bdX$ is disconnected then we still assume a fibration $\phi:\bdX\to B$ is given, but $B$ may now be disconnected.\footnote{The fibres can be disconnected even if $\bdX$ is connected, e.g.\ if $\phi:S^1\to S^1$ is the double cover.}
Note that fibres over different components of $B$ need not be diffeomorphic, and these components can even have different dimensions.
In addition, for $c-\phi$-metrics the value of $c$ may vary from component to component. In results that allow for different values of $c$ (e.g.\ in Section \ref{sec:calderon}) this means that they apply to spaces having singularities (or ends) of different types, like an incomplete cusp singularity as well as an asymptotically conical end.
\end{remark}

\subsection{Examples of fibred cusp spaces}
\label{subsec:examples}

Fibred cusp spaces appear in many different contexts. We give a number of examples. For easy reading we indicate the one-forms from which the metrics are built as positive definite quadratic forms, locally near $x=0$, and leave out indices. These one-forms are $x^c$ times the one-forms \eqref{eqn:phi-forms}. Note that sometimes the $z$-coordinates are missing, which means that the fibre is a point, and sometimes the $y$-coordinates are missing, which means that the base is a point. In some cases conditions (like Riemannian submersion etc.) should be imposed.

Note that $c-\phi$-metrics are complete if $c\leq1$ and incomplete if $c>1$.

\smallskip

\noindent\boxed{$c=0$:}

\begin{description}
[align=left]
\item
[$\displaystyle\frac{dx}{x^2}, \frac{dy}x$] asymptotically conic (aka scattering) metrics; see Figure \ref{fig:examples} right\footnote{Various subclasses of these have been studied, like asymptotically Euclidean (AE) or asymptotically locally Euclidean (ALE) metrics.  These correspond to $B=S^{n - 1}$ with the round metric -- the simplest example being $X=\overline{\R^n}$, the radial compactification of $\R^n$ with the Euclidean metric --  and $B=S^{n-1}/\Gamma$ for a finite subgroup $\Gamma\subset O(n,\mathbb{R})$ acting freely, respectively. They occur, for example, in general relativity, as solitons in Ricci flow and as gravitational instantons (complete hyper-K\"ahler 4-manifolds \cite{kronheimer}).} 
\\
\item
[$\displaystyle\frac{dx}{x^2}, dz$] infinite cylinder $(1,\infty)_r\times F$, compactified at infinity to $(1,\infty]\times F$ with smooth structure on $(1,\infty]$ defined using $x=\frac1r$ \footnote{That is, a function $u:(1,\infty)\to\C$ is defined to be smooth on $(1,\infty]$ if $u(r)=v(\frac1r)$ where $v:(0,1)\to\C$ extends smoothly to $[0,1)$.}
\item
[$\displaystyle\frac{dx}{x^2}, \frac{dy}x, dz$] a simple example is $\overline{\R^n}\times F$ for $(F,g_F)$ a compact Riemannian manifold; other classes with $F=S^1$ are given by asymptotically locally flat (ALF) gravitational instantons (where $B=S^2/\Gamma$) and asymptotically locally conical (ALC) metrics with $G_2$-holonomy, where $b=5$ \cite{FosHasNor}
\end{description}

\smallskip

\noindent\boxed{$c=1$:}

\begin{description}
[align=left]
\item[$\displaystyle\frac{dx}x, dy$]  infinite cylinder $(1,\infty)_r\times F$, compactified at infinity to $(1,\infty]\times F$ with smooth structure defined using $x=e^{-r}$
\item[$\displaystyle\frac{dx}x, x\,dz$] hyperbolic cusp, see Figure \ref{fig:examples} center \\
\item [$\displaystyle\frac{dx}x, dy, x\,dz$] certain locally symmetric spaces (see \cite{Bor:SRCAG}, \cite{Mue:MCRO} and the generalizations below);  K\"ahler-Einstein metrics on certain quasiprojective varieties, see \cite{RocZha:ACKMFVQM}

\end{description}

\smallskip

\noindent\boxed{$c=2$:}

\begin{description}
[align=left]
\item [$dx, x\,dy$] cone near the cone tip
\item [$dx, x^2\,dz$] incomplete cuspidal singularity, see Figure \ref{fig:examples} left\footnote{\label{fn:other cusps}%
Note that if we allow non-integer $c>1$ then replacing $x$ by $\xbar = x^{c-1}$ (of course this changes the smooth structure at $x=0$) then $x^c\frac{dx}{x^2}$, $x^c\, dz$ turn into
$d\xbar, \xbar^\gamma \,dz $
(up to constant factors) where $\gamma = \frac c{c-1}$, so we obtain 'sharper' incomplete cusps (i.e., $\gamma>2$) if $1<c<2$ and 'blunter' ones (i.e., $1<\gamma<2$) if $c>2$.
Cusps with any 'sharpness' exponent $\gamma>1$ arise, for example, as the subset $\{(u,v,\xbar): u^2+v^2 = \xbar^{2\gamma},\xbar\geq0\}\subset\R^3$  (up to mixed terms).
}
\item [$dx, x\,dy, x^2\,dz$] a simple example (where, however, $F$ has a boundary, see Section \ref{sec:calderon}) is the complement of two disjoint balls in $\R^3$, touching at a point, see Figure \ref{fig:blow-up example}

\end{description}
\begin{remark}
\label{rem:meaning phi}
Given this large number of example classes the reader might wonder which spaces (with sufficiently regular structure near infinity) are \textit{not} fibred cusp spaces. In order to answer this, note that the common feature of $c-\phi$-metrics over all different values of $c$, preserved under changes of the $x$-variable like $\xbar=x^\alpha$, $\alpha>0$, \footnote{Here we may also allow non-integer $\alpha$ and $c$ modulo minor smoothness issues, compare footnote \ref{fn:other cusps}.} is that the total exponent of $x$ in the $dx$ term, which is $x^{c-1}$, is less than or equal to the exponents of $x$ in all other terms.\footnote{Really we are talking about the scaling exponent in $x$: This is $c-1$ 
for the one-form $x^c\frac{dx}{x^{2}}$ since replacing $x$ by $\lambda x$ multiplies it by $\lambda^{c-1}$.}
\footnote{This remains true even for the coordinate change $\xbar=e^{-1/x}$, which turns $\frac{dx}{x^2}$ into $\frac{d\xbar}{\xbar}$
(compare the second example above for $c=0$ with the first example for $c=1$), in the sense that the coefficients of the $dy$, $dz$ terms are smaller, as $x\to0$, than the total $x$-power of the $dx$ term.}
Geometrically, this can be expressed via the induced one-forms $x^{c-1}dy$, $x^c dz$ on the hypersurfaces $H_{x_0}=\{x=x_0\}$ near the boundary, which determine lengths (e.g., diameter) and volume form on $H_{x_0}$: Using that $\int x^c\,\frac{dx}{x^2}$ equals $\const\, x^{c-1}$ if $c\neq1$ and $\log x$ if $c=1$ we get (always in the regime $x\to0$):
\begin{itemize}
\item  $c<1$ (\lq very complete case\rq): Here $H_x$ is roughly the sphere of large radius $r=\const\cdot x^{c-1}$ around any interior point $p\in\interior{X}$. The induced 1-forms are $r\,dy$ and $r^\gamma\,dz$ where $\gamma=\frac c{c-1}=1-\frac1{1-c}<1$.
So the length element on large spheres grows at most like the radius $r$ as $r\to\infty$, and more slowly in fibre directions.
\item $c=1$ (\lq marginally complete case\rq): Similar to $c>1$, but the radius is $r\approx \log\frac1x$ and the one-forms are $dy$ and $e^{-r}dz$: No growth of the base, exponential collapse of the fibres.
\item $c>1$ (incomplete case): Now geometrically the boundary is collapsed to a point\footnote{This means that any two points of $\bdX$ have distance zero to each other.}, and  $H_x$ is roughly the sphere of \textit{small} radius $r=\const\cdot x^{c-1}$ around this point. The induced 1-forms are $r\,dy$ and $r^\gamma\,dz$ where now $\gamma=\frac c{c-1}>1$. So the spheres shrink linearly, with fibres collapsing faster. See Figure \ref{fig:blow-up example} for an example.
\end{itemize}
Summarizing, \textbf{the essence of fibred cusp spaces} is, roughly\footnote{This extends to the more general multiply fibred cusp setting, see \eqref{eqn:mult fibred phi}, and sets it apart from the edge setting with non-trivial base, see Footnote \ref{fn:edge}. Here we are neglecting, for instance, smoothness issues at the singularity/at infinity.}:
In the complete case the linear size of large spheres grows at most linearly in the radius; in the incomplete case the linear size of small spheres centered at the singularity shrinks at least linearly.
\end{remark}

\subsection{Generalizations: Higher order, iterated and multiply fibred cusps}

Our setting has various natural generalizations: 
\begin{itemize} 
\item higher order and multiply fibred cusps: 
locally we can replace \eqref{eqn:phi-forms} by 
\begin{equation}
\label{eqn:mult fibred phi}
\frac{dx}{x^{b_0+1}}\,,\ \frac{dy_0}{x^{b_0}}\,,\ \frac{dy_1}{x^{b_1}}\,,\dots,\ dy_m\,,
\end{equation}
with each $y_i$ a group of variables and integers $b_0>\dots>b_{m-1}>0$; as before, we may also consider these one-forms multiplied by the conformal factor $x^{c}$.
The geometric data would be a stack of fibrations
$\partial X=B_m\to B_{m-1}\to\dots\to B_0$ (with $y_0\in B_0$, $y_1$ in the fibre of $B_1\to B_0$ etc.). See \cite{GriHun:POCGQLSSI}. For $m=2$ and $c=b_0$ such metrics arise on locally symmetric spaces of $\Q$-rank 1, where $b_1$ and $b_2$ are related to Lie-theoretic weights, see \cite{Bor:SRCAG}, \cite{Mue:MCRO}. Remark \ref{rem:meaning phi} applies accordingly, where the $b_i$ may be real.
\item iterated (multiply) fibred cusps: here $X$ is a manifold with corners, and  each boundary hypersurface is equipped with a fibration of manifolds with corners (or a stack of fibrations), and certain compatibility conditions are imposed at the corners. This structure is called a manifold with fibred corners. \footnote{
For metrics with $c=1$ see \cite{DebLesRoc:POMFC}. 
For $c=0$ iterated versions of many classes of $\phi$-metrics commonly get the modifier \lq  quasi\rq. For instance, QALE metrics (quasi-asymptotically locally Euclidean, i.e. asymptotic to $\R^n/\Gamma$ with $\Gamma$ acting not freely on $S^{n-1}$ and behaving in prescribed ways near the singularities) generalize ALE metrics, see e.g.\ the instantons in \cite{joyce}. Scattering (aka AC) metrics (which include ALEs) are generalized by QAC (quasi-asymptotically conical) metrics, with certain non-compact Calabi-Yau manifolds \cite{ConDegRoc} and Nakajima quiver varieties \cite{DimRoc} as examples (and a \lq warped\rq\ variant in   \cite{ConRoc}). 
All of these have point fibres. With general fibres one gets QFB (quasi-fibred boundary, see \cite{KotRoc:QFBPO}) metrics which generalize $0$-$\phi$-metrics; they appear in the compactification of the moduli spaces of $\text{SU}(2)$-monopoles of magnetic charge $3$ (those of charge $2$ are ALF) and are expected to model the $L^2$-metric in the moduli space of $\text{SL}(2,\mathbb{C})$ Higgs bundles, up to some singular fibres \cite{KotRoc}. They were thus partly introduced to solve (parts of) the Sen and Vafa-Witten conjectures.}
\footnote{Other examples, for certain positive $c$, are locally symmetric spaces of higher $\Q$-rank. Here $X$ is the Borel-Serre compactification \cite{BorJi:CSLSS}.}
\item
edge spaces: these are manifolds with fibred boundary, but \eqref{eqn:phi-forms} is replaced by $\frac{dx}{x}, \frac{dy_i}x, dz_j$ (or more generally by \eqref{eqn:mult fibred phi} with the first term replaced by $\frac{dx}{x^{b_0}}$ -- this general setting would include fibred cusp spaces as the special case where $B_0$ is a point); for $c$-edge spaces we multiply by the conformal factor $x^c$. 1-edge metrics are also called wedge metrics.
Examples where $F$ is a point are hyperbolic funnels ($c=0$) and compact manifolds with boundary and smooth metric ($c=1$). 
\footnote{\label{fn:edge}%
There is some overlap between the classes of $c-\phi$-metrics and $c'$-edge-metrics, for example cylindrical ends are also $0$-edge with trivial base.
However, $c'$-edge metrics with non-trivial base are fundamentally different from all $c-\phi$-metrics: the base for the former is much larger as $x\to0$ since $dy$ has an extra factor of $x^{-1}$ in comparison. More precisely, in reference to Remark \ref{rem:meaning phi}, in the very complete case, $c'<0$, the radius is now $r\approx x^{c'}$, and the coefficient of $dy$ is $r^{\beta}$ where $\beta=\frac{c'-1}{c'}>1$, rather than $r$. In the marginally complete case $c'=0$, the radius is $r\approx \log\frac1x$, and the coefficient of $dy$ is $e^r$, which corresponds to exponential volume growth of spheres in hyperbolic space. In the incomplete case, $c'>0$, we have again $r=x^{c'}$, which is now small, and the coefficient of $dy$ is $r^{\beta}$ where $\beta=1-\frac1{c'}<1$, so balls shrink more slowly than $r$ as $r\to0$ (or not at all if $c'=1$, and even blow up if $0<c'<1$).
}

The theory for these structures is quite different from the theory for fibred cusp spaces. See \cite{mazzeo:edge1} \cite{mazzeo-vertman:edge2}. In \cite{albgr:pseudo}, an analysis of the spectral and index theory of iterated wedges is carried out.

\end{itemize}

\section{The $\phi$-calculus}
\label{sec:phi-calculus}
In Section \ref{sec:fib cusp} we defined $\phi$-manifolds, $\phi$-vector fields, and (general) $\phi$-metrics. We now introduce $\phi$-differential operators and the $\phi$-pseudodifferential calculus, and then illustrate the concepts by the simple example of the Laplacian on $\R^n$.
The natural differential operators associated to a general $\phi$-metric  (e.g. Laplacian, Dirac operator) are of a special type, called \textbf{$\phi$-differential operators}. By definition, these are obtained by composing $\phi$-vector fields and adding functions. So they have the form, in coordinates near $\bdX$,
\begin{equation}
 \label{eqn:phi operator}
 P = \sum_{k+|\alpha|+|\beta|\leq m} a_{k,\alpha,\beta}(x,y,z) (x^2D_x)^k (xD_y)^\alpha D_z^\beta,\quad a_{k,\alpha,\beta} \text{ smooth, }
\end{equation}
where $D_x=\frac1i\partial_x$ etc., $\alpha,\beta$ are multi-indices and $m\in\N_0$ is the order of $P$. Importantly, smoothness is assumed up to the boundary $x=0$. The set of these operators is denoted $\Diff_\phi^m(X)$.  The class of operators acting between sections of vector bundles $E,E^\prime$ on $X$ is defined similarly and denoted $\Diff^m_\phi(X;E,E^\prime)$. 
Everything that follows extends to that case in a straightforward way, but to simplify notation we stick to the scalar case.
The $\phi$-\textbf{principal symbol} of $P\in\Diff^m_\phi(X)$ is the standard principal symbol in the interior, and near $\bdX$ where \eqref{eqn:phi operator} holds it is the function
$$\phisigma_m(P) =  \sum_{k+|\alpha|+|\beta| = m} a_{k,\alpha,\beta}(x,y,z) \tau^k \eta^\alpha \zeta^\beta$$
for $\tau\in\R$, $\eta\in\R^b$, $\zeta\in\R^f$.
Invariantly and globally, $\phisigma_m(P)$ can be considered as a function  (or section of a suitable vector bundle of homomorphisms) on  $\phiT^*X$. The operator $P$ is called $\phi$-\textbf{elliptic} if its $\phi$-principal symbol is invertible outside the zero section $(\tau,\eta,\zeta)=0$.

The \textbf{normal family} of $P\in \Diff^m_\phi(X)$ captures its behavior at the  boundary, globally in the fibres. It is a family of differential operators
on the fibres $F_y=\phi^{-1}(y)$, parametrized by $y\in B$, $\tau\in\R$ and $\eta\in T_y^*B$. In local coordinates where \eqref{eqn:phi operator} holds it is given by\footnote{In the literature the notation $\hat N(P)$ is sometimes used for the normal family, with $N(P)$ denoting the normal \emph{operator}, where $\tau,\eta$ are replaced by differentiations $D_T$, $D_Y$ in auxiliary variables $T\in\R$, $Y\in T_yB$. We do not use the normal operator here.}
\begin{equation}
\label{eqn:def normal op diff}
 N(P)(\tau;y,\eta) =  \sum_{k+|\alpha|+|\beta|\leq m} a_{k,\alpha,\beta}(0,y,z) \tau^k \eta^\alpha D_z^\beta\,.
\end{equation}
A $\phi$-differential operator $P$ is called \textbf{fully elliptic} if it is $\phi$-elliptic and $N(P)(\tau;y,\eta)$ is invertible for all $\tau,y,\eta$ (including $\tau=0$, $\eta=0$). 

Note that the Laplacian for a $c-\phi$-metric is of the form $x^{-2c}P$ where $P\in\Diff^2_\phi(X)$ is $\phi$-elliptic, and similarly for the Dirac operator.

\medskip

We now  describe the $\phi$-pseudodifferential calculus introduced in \cite{MM98}. This is a pseudodifferential calculus, i.e.\ a set of operators closed under composition, which contains $\Diff_\phi^*(X)$ as well as parametrices of fully elliptic elements of this space. By the general philosophy on singular pseudodifferential calculi introduced by R. Melrose,
such a calculus is defined by a set of Schwartz kernels, which are distributions on the (interior of the) double space $X^2=X\times X$, and whose boundary behavior is restricted by requiring that their pull-backs to a suitable blow-up, $X^2_\phi$, of $X^2$ satisfy certain smoothness and vanishing conditions at the boundary hypersurfaces of $X^2_\phi$ and have a conormal singularity at the diagonal, uniformly up to the boundary.\footnote{See also  \cite{KotRoc:PMWFC} and \cite{Hin:POMWSBG} for new general perspectives on such pseudodifferential calculi.
Yet another perspective is via Lie groupoids; the relation between both formalisms for \textit{iterated} fibred cusp metrics (of which our setting is the simplest instance -- depth 1) is analysed in \cite{DebLesRoc:POMFC}.
The Lie groupoid integrating the Lie algebroid $\leftindex^\phi TX$ is the $\phi$-double space $X^2_\phi$ with its lateral faces removed.}

The $\phi$-double space, $X^2_\phi$, is defined as follows, see Figure \ref{fig:phi double space}. We first assume $\bdX$ is connected (but see Remark \ref{rem:connected2}).
 First, blow up $(\dX)^2\subset X^2$. This yields the b-double space with blow-down map
\begin{equation}
\label{eqn:b-blowup}
\beta_b: X^2_b := [X^2;(\dX)^2] \to X^2 \,.
\end{equation}
Its front face is denoted by $\bbf$, and the interior of $\bbf$ is naturally diffeomorphic to $(0,\infty)\times(\dX)^2$ where the first coordinate is $t=\frac{x'}x$, with $x,x'$ the pull-backs of the boundary defining function $x$ on $X$ to the first and second factor in $X^2$.
Let $\diag_X\subset X^2$ be the diagonal and $\diag_{b}$ be its lift to $X^2_b$. It meets the boundary of $X^2_b$ in the interior of $\bbf$, in the set $\{1\}\times\diag_\dX$. The larger interior submanifold $D_\phi:=\{1\}\times\{(p,p')\in(\dX)^2:\phi(p)=\phi(p')\}$ of $\bbf$ is called the fibre diagonal.
It can also be described as the intersection of the lift of the interior fibre diagonal,
\begin{align}
\label{def diag phi int}
&\textup{diag}_{\phi,\mathrm{int}}= \{ (p, p') \in {\cU} \times 
{\cU}: \phi(p) = \phi(p')\},
 \end{align}
with $\textup{bf}$; that is, $D_\phi=\beta_b^*(\textup{diag}_{\phi,\mathrm{int}}) \cap \textup{bf}$.

We blow  up $D_\phi$ and define the $\phi$-double space
\begin{equation}
\label{eqn:phi-blowup}
\beta_{D_\phi}: X^2_\phi := [X^2_b;D_\phi]\to X^2_b,\quad  \beta_\phi:= \beta_b\circ\beta_{D_\phi}:X^2_\phi\to X^2\,.
\end{equation} 
The front face created by this blow-up is denoted $\phif$. The lift of $\bbf$ is still denoted $\bbf$.
The diagonal $\diag_{b}$ lifts to a p-submanifold $\diag_\phi$ of $X^2_\phi$ which meets the boundary in the interior of $\phif$.

We describe these spaces locally, using  local coordinates $x,y,z$. Starting from coordinates $(x,y,z; x',y',z')$ on $X^2$
 we have $(\dX)^2= \{x=x'=0\}$, so coordinates near interior points of $\bbf$ are
 $$ t=\frac {x'}{x}\,,\ x\,,\ y\,,\ z\,,\ y'\,,\ z'\,,$$
 with $x$ defining $\bbf$ there.
The diagonal in $X^2$ is $\{x=x',y=y',z=z'\}$ and lifts to $\{t=1,y=y',z=z'\}$.
 The fibre diagonal $D_\phi$ is locally $\{t=1, x=0, y=y'\}$. So projective coordinates near the interior of $\phif$ are
\begin{equation}
\label{eqn:phi double coords}
 T := \frac{t-1}{x}\,,\ x\,,\ Y:=\frac{y-y'}{x}\,,\ y\,,\ z\,,\ z',
\end{equation}
 with $x$ defining $\phif$ there.
 In these coordinates the diagonal $\diag_\phi$ is given by $\{T=0, Y=0, z=z'\}$, and $|(T,Y)|\to\infty$ corresponds to the boundary of $\phif$, which is its intersection with the lift of $\bbf$.
\begin{remark}[Geometric meaning of $X^2_\phi$]
\label{rem:X2phi geometry}
There are two interrelated ways to understand the significance of the blow-up $\beta_\phi:X^2_\phi\to X^2$ geometrically:
\begin{enumerate}
\item
The $\phi$-vector fields on $X$, when considered as vector fields on $X^2$ in the $x,y,z$ variables, i.e.\ pulled back from the left factor, lift under $\beta_\phi$ to smooth vector fields on $X^2_\phi$ that span a vector bundle of rank $\dim X$, whose fibre at any $\gamma\in \diag_\phi$ is transversal to $T_\gamma\diag_\phi$ as 
 subspace of $T_\gamma X^2_\phi$.
\item
The distance function $d:\interior X\times\interior X\to\Rplus$ for any $\phi$-metric is resolved by $\beta_\phi$, in the following sense. Denote $\dtilde=\beta_\phi^*d$.
\begin{enumerate}
\item[(a)]
For $p,p' \in \interior X$ 
approaching $\dX$ the distance $d(p,p')$ remains bounded iff $\beta_\phi^{-1}(p,p')$ approaches $\phif$ in its interior.

More precisely, the function $1/\dtilde$ extends smoothly to $X^2_\phi\setminus\diag_\phi$, and the extension vanishes precisely at $\partial(X^2_\phi)\setminus\interior{\phif}$.
  
\item[(b)] 
The function ${\tilde d}\,{}^2$ extends smoothly to $(X^2_\phi)^\circ \cup \interior{\phif}$, vanishing quadratically at $\diag_\phi$ with positive definite transversal Hessian, uniformly to $\diag_\phi\cap\;\phif$.
\end{enumerate}
Similar statements hold for the distance function $d_c$ associated with a $c$-$\phi$-metric, if one replaces $\dtilde_c=\beta_\phi^*d_c$ by $x^{-c}\dtilde_c$ (assuming $x\leq 2x'$, i.e.\ in the upper two-thirds of the picture in Figure \ref{fig:phi double space}, and similarly with the roles of $x,x'$ interchanged).\footnote{More precisely, we have the following leading orders of $\dtilde_c$ at the boundary hypersurfaces\\
 $\phif$, $\bbf$, $\lf$ (left face, black arrow in Figure \ref{fig:phi double space}), with non-vanishing coefficient in their interiors:\\ 
If $c=0$: $x^0, x^{-1}, x^{-1}$. If $c=1$: $x, x^0, \log x$. If $c=2$: $x^2, x^1, x^0$.
The logarithmic divergence at $\lf$ for $c=1$ corresponds to the typical logarithmic length scale in hyperbolic space.
}
In particular, for $p,p'$ approaching $\partial X$
we have $d_c(p,p') = O(x^c)$ iff  
$\beta_\phi^{-1}(p,p')$ approaches $\phif$ in its interior.
\end{enumerate}
\end{remark}
\begin{figure}[!ht]
\centering
\resizebox{1\textwidth}{!}{%
\begin{circuitikz}
\tikzstyle{every node}=[font=\LARGE]
\draw [->, >=Stealth] (2.5,17.5) -- (2.5,25);
\draw [->, >=Stealth] (2.5,17.5) -- (10,17.5);
\draw [ color={rgb,255:red,228; green,7; blue,7}, line width=1pt, dashed] (2.5,17.5) -- (10,25);
\draw [ color={rgb,255:red,247; green,145; blue,2}, line width=1pt, short] (2.5,17.5) .. controls (3.75,19) and (7.75,22.5) .. (7.5,25);
\draw [ color={rgb,255:red,19; green,111; blue,231}, line width=1pt, short] (2.5,17.5) .. controls (4.5,19.55) and (7,21.25) .. (10,21.25);
\draw [ color={rgb,255:red,21; green,157; blue,12}, line width=1pt, short] (2.5,17.5) -- (5,25);
\draw [ color={rgb,255:red,115; green,11; blue,218}, line width=1pt, ->, >=Stealth] (3.75,23.75) -- (2.5,23.75);
\node [font=\LARGE] at (10,17) {\fontsize{30}{28}\selectfont $x$};
\node [font=\LARGE] at (2,24.75) {\fontsize{30}{28}\selectfont $x'$};
\draw [->, >=Stealth] (15,21.25) -- (15,25);
\draw [->, >=Stealth] (18.75,17.5) -- (22.5,17.5);
\draw [short] (18.75,17.5) .. controls (18.75,19.75) and (17.25,21.25) .. (15,21.25);
\draw [ color={rgb,255:red,115; green,11; blue,218}, line width=1pt, ->, >=Stealth] (16.25,23.75) -- (15,23.75);
\draw [ color={rgb,255:red,21; green,157; blue,12}, line width=1pt, short] (16.25,21.1) -- (17.5,25);
\draw [ color={rgb,255:red,228; green,7; blue,7}, line width=1pt, dashed] (17.6,20.35) -- (22.5,25);
\draw [ color={rgb,255:red,239; green,144; blue,0}, line width=1pt, short] (17.6,20.35) -- (20,25);
\draw [ color={rgb,255:red,42; green,113; blue,233}, line width=1pt, short] (17.6,20.35) -- (22.5,21.25);
\node at (2.5,17.5) [circ, color={rgb,255:red,223; green,12; blue,184}] {};
\node at (17.6,20.35) [circ, color={rgb,255:red,213; green,18; blue,185}] {};
\draw [->, >=Stealth] (13.75,21.25) -- (11.25,21.25);
\draw [->, >=Stealth] (26.25,21.25) -- (23.75,21.25);
\draw [->, >=Stealth] (27.5,21.25) -- (27.5,25);
\draw [->, >=Stealth] (31.25,17.5) -- (35,17.5);
\draw [ color={rgb,255:red,115; green,11; blue,218}, line width=1pt, ->, >=Stealth] (28.75,23.75) -- (27.5,23.75);
\draw [ color={rgb,255:red,21; green,157; blue,12}, line width=1pt, short] (28.75,21.1) -- (30,25);
\draw [ color={rgb,255:red,228; green,7; blue,7}, line width=1pt, dashed] (30.5,20.5) -- (35,25);
\draw [ color={rgb,255:red,239; green,144; blue,0}, line width=1pt, short] (30.25,20.65) -- (32.5,25);
\draw [ color={rgb,255:red,42; green,113; blue,233}, line width=1pt, short] (30.75,20) -- (35,21.25);
\draw [short] (31.25,17.5) .. controls (31.25,18.25) and (31,18.75) .. (30.75,19.25);
\draw [short] (27.5,21.25) .. controls (28.25,21.25) and (28.75,21.25) .. (29.25,20.75);
\draw [short] (30.75,19.25) .. controls (31,20.25) and (30.25,21) .. (29.25,20.75);
\draw [->, >=Stealth] (18,18.75) .. controls (18,19.5) and (16.75,20.5) .. (16.25,20.5) ;
\node [font=\LARGE] at (17,19.5) {\fontsize{30}{28}\selectfont $t$};
\draw [->, >=Stealth] (30.5,19.5) .. controls (30.5,20.5) and (29.8,20.55) .. (29.5,20.5) ;
\node [font=\LARGE] at (29.75,19.75) {\fontsize{30}{28}\selectfont $T$};
\node [font=\Huge] at (6.5,15.25) {\fontsize{30}{28}\selectfont $X^2$};
\node [font=\Huge] at (19,15.25) {\fontsize{30}{28}\selectfont $X^2_b$};
\node [font=\Huge] at (31.5,15.25) {\fontsize{30}{28}\selectfont $X^2_\phi$};
\end{circuitikz}
}%

\caption{$\phi$-double space depicted for $X=\mathbb{R}_+$ (so $\partial X =\{0\}$ has trivial base and fibre). First the corner $\{(0,0)\}=\left(\partial X\right)^2\subset X^2$ is blown up, which separates lines arriving at the origin from different directions (note how the green is separated from the orange and blue ones by the new face with coordinate $t=\frac{x'}{x}$, since curves $x'=t_1x + O(x^2)$ and $x' = t_2x + O(x^2)$ land at the point of the front face with $t=t_1$ or $t=t_2$, respectively). Moreover, different parabolae arriving at the origin from the same direction (like orange and blue) land in the same point at the front face but look straightened up in the neigborhood of it. In particular, if they approach the origin tangent to the diagonal $x' = x$ (red, dashed), they will land at $t=1$. Blowing up this point of incidence $\{t=1, x=0\}\subset X^2_b$ separates the different directions of approach to it, so in the new front face appearing in $X^2_\phi$, the lifts of the orange, red and blue lines are now separated. This corresponds to introducing a coordinate $T=\frac{t-1}x=\frac{x'-x}{x^2}$ that distinguishes the second order behaviour, i.e. curves $x'=x+Tx^2+ O(x^3)$ with different $T$'s.\\ In terms of operators, an element $P\in\Psi^m_\phi(\mathbb{R}_+)$ has a kernel $K_P$ that is $O(x^\infty)$ as $x\to 0$ while $x'>\varepsilon$ (purple) or  as $x\to 0$ and $\frac{x'}{x}=t +O(x)$ with $t\neq 1$ (green), and is smooth as a function of $T\neq0$ and $x\geq0$ when setting $x'=x + Tx^2$ (e.g. orange and blue), and has a conormal singularity at $T=0$ (red), uniformly to $x=0$.
\\
If $B$, $F$ are non-trivial, with variables $y$, $z$, then there are additional $y,z,y',z'$ variables (transversal to the drawing plane), and the $\phi$-blow-up is only done at $y=y'$.
}
\label{fig:phi double space}
\end{figure}

We consider operators $P$ acting on functions on $\interior{X}=X\setminus\dX$ which are given by Schwartz kernels $K_P$, which are distributions on $(\interior X)^2=\interior{X^2}$, in the sense that
\begin{equation}
 \label{eqn:schwartz kernel}
 (Pu)(p) = \int_{\interior X} K_P(p,p') u(p')\, \nu(p'),
\end{equation}
where $\nu$ is some fixed density on $\interior X$. We choose (and fix once and for all) for $\nu$ a smooth positive $\phi$-density, i.e.\ locally $\nu = a \frac{dx}{x^2}\frac{dy}{x^b}dz$ with $a>0$ smooth up to the boundary $\{x=0\}$ and $b = \dim B$.
For example, $\nu$ could be the volume density of a $\phi$-metric. The reason for this choice is that for $P=\Id$ (and, say, $a=1$), we have $K_P=(x')^{b+2}\delta(x-x')\delta(y-y')\delta(z-z')$, which in coordinates \eqref{eqn:phi double coords} is $\delta(T)\delta(Y)\delta(z-z')$, and since this has no $x$-factor, it extends from the interior of $X^2_\phi$  to a distribution on $X^2_\phi$ as a smooth non-vanishing delta-type distribution for the submanifold $\diag_\phi$.
The lifting property of $\phi$-vector fields mentioned in Remark \ref{rem:X2phi geometry}  then implies that the Schwartz kernels for $P\in\Diff^m_\phi(X)$ lift to $\interior{X^2_\phi}$ and extend to $X^2_\phi$ to be delta-type distributions of order at most $m$ for $\diag_\phi$.
Explicitly, if $P$ is given in coordinates by \eqref{eqn:phi operator} and $\nu=\frac{dx}{x^2}\frac{dy}{x^b}dz$ then
\begin{equation}
 \label{eqn:kernel of phidiff op}
 K_{P} = \sum_{k+|\alpha|+|\beta|\leq m} a_{k,\alpha,\beta}(0,y,z) D_T^k\delta(T) D_Y^\alpha \delta(Y) D_z^\beta\delta(z-z') + O(x).
\end{equation}
We define $\phi$-pseudodifferential operators by replacing delta-type distributions by the larger space of classical conormal distributions, defined in the Appendix:

\begin{definition}
Let $X$ be a $\phi$-manifold and $m\in\R$. The space $\Psi_\phi^m(X)$ of $\phi$-pseudodifferential operators of order $m$ is defined as the set of operators \eqref{eqn:schwartz kernel} whose Schwartz kernels $K_P$ lift to $X^2_\phi$ to elements of
$$ I^m_\cl(X^2_\phi, \diag_\phi) $$
that vanish to infinite order at all boundary hypersurfaces of $X^2_\phi$ except $\phif$.\footnote{The special case where the fibre is a point is sometimes called scattering calculus and denoted $\Psi^m_{\text{sc}}(X)$, see \cite{Mel:SSTLAES}, \cite{Mel:GST}; it is closely related to the so-called SG-calculus, see \cite{CorRod:ESGC}. The special case where the base is a point is sometimes called cusp calculus (see \cite{SchSteSha:OAMWCS} for a different approach to this setting).}
\end{definition}
See Figure \ref{fig:phi double space}.
There is also a more general definition without $\cl$, but all operators occurring in this paper are classical. By the remarks before the definition we have $\Diff_\phi^*(X) \subset \Psi_\phi^*(X)$. The definition of the $\phi$-principal symbol extends in a straight-forward way using the local representation \eqref{eqn:def conormal}. The definition of the normal family of $P\in\Diff_\phi^m(X)$ given in \eqref{eqn:def normal op diff} does not extend directly to $\Psi_\phi^*(X)$. However, reinterpreting this formula in terms of \lq oscillatory testing\rq\ allows to extend it as follows: for $\tau\in\R$, $\eta\in\R^b$ and $y_0\in B$ let $g(x,y)=-\frac\tau x + \frac\eta x (y-y_0)$ in coordinates on $B$ near $y_0$.
Then $x^2D_x e^{ig} = \left(\tau-\eta(y-y_0)\right)e^{ig}$ and $xD_y e^{ig} = \eta e^{ig}$, and this implies that for $P$ as in \eqref{eqn:phi operator} and $u\in\Cinf(X)$ we have
\begin{equation}
 \label{eqn:normal op osc testing}
 \left[e^{-ig}P(e^{ig}u)\right]_{|F_{y_0}} = N(P)(\tau;y_0,\eta)(u_{|F_{y_0}})
\end{equation}
as functions on the fibre $F_{y_0}$, since $F_{y_0}$ is given by $x=0, y=y_0$.
It can be shown that the left hand side is well-defined and smooth for $P\in \Psi^*_\phi(X)$ and only depends on $u_{|F_{y_0}}$, and that $N(P)$ so defined is a (standard) pseudodifferential operator with parameter $(\tau,
\eta)$ on $F_{y_0}$, and varies smoothly in $y_0$. Also, the definition makes sense invariantly when considering $\eta\in T_{y_0}B$. By \eqref{eqn:def normal op diff} the normal family for $P\in\Diff^m_\phi(X)$ vanishes if and only if $P\in x\Diff^m_\phi(X)$, and an analogous statement holds for $P\in\Psi^* _\phi(X)$.
A short calculation shows that the Schwartz kernel $K_{N(P)}(\tau;y,\eta;z,z')$ of $N(P)(\tau;y,\eta)$ is the $(T,Y)\to(\tau,\eta)$ Fourier transform of the restriction of $K_P$ to $\phif$, when writing $K_P$ in coordinates \eqref{eqn:phi double coords}.

As for $\phi$-differential operators, a $\phi$-pseudodifferential operator is called $\phi$-elliptic if its $\phi$-principal symbol is invertible outside the zero section and fully elliptic if in addition its normal family is invertible for all $\tau,y,\eta$.

We denote by $L^2_\phi(X)$ the $L^2$-space with respect to the volume form of a $\phi$-metric, and define $\phi$-Sobolev spaces $H^k_\phi(X):=\{u:\,\Diff^k_\phi(X)u\subset L^2_\phi(X)\}$ for $k\in\N_0$. Sobolev spaces $H^s_\phi(X)$ are then defined for $s\in\R$ by interpolation and duality.

The main facts about the $\phi$-pseudodifferential calculus are:
\begin{enumerate}
\item $\Psi^*_\phi(X):= \bigoplus_{m\in\R}\Psi^m_\phi(X)$ is an $\R$-graded $\star$-algebra, i.e.\ a vector space and closed under adjoints and composition, with orders adding under composition.
\item The $\phi$-principal symbol $P\mapsto \phisigma(P)$ and the normal family $P\mapsto N(P)$ are $\star$-algebra homomorphisms%
\footnote{The $\phi$-principal symbol of $P\in \Psi^m_\phi(X)$ lies in $S^{[m]}({}^\phi \dot{T}^*X)$, the space of smooth functions on ${}^\phi \dot{T}^*X = {}^\phi {T}^*X\setminus\{\text{zero section}\}$ which are positively homogeneous of degree $m$ in the fibres. The normal family lies in the algebra of families of pseudodifferential operators on the fibres $F_y$ of $\phi$ with parameter $(\tau,\eta)\in {}^\phi NB$, and depending smoothly on $y\in B$. See \cite{MM98}.}, i.e.\ they are linear and respect composition, and $ \phisigma(P^*) = \overline{ \phisigma(P)}$ where $P^*$ is the adjoint with respect to any $\phi$-metric.
\item Operators in $\Psi^m_\phi(X)$ are bounded $H^s_\phi(X)\to H^{s-m}_\phi(X)$ for all $s$ and map each of the spaces $\calA^a(X)$, $\calA^{\calE}(X)$ (cf. Appendix) into itself, for any $a\in\R$ and index set $\calE$.
\item An element of $x^a\Psi^m_\phi(X)$ is a compact operator in $L^2_\phi(X)$ if and only if $m<0$ and $a>0$.
\item An operator $P\in \Psi^m_\phi(X)$ is fully elliptic if and only if it has a parametrix with remainders in $x^\infty\Psi^{-\infty}_\phi(X)$, if and only if it is Fredholm as a map $H^s_\phi(X)\to H^{s-m}_\phi(X)$ for any $s$.
In particular, $\ker P\subset \Cdotinf(X)$ in this case (i.e., its elements vanish to infinity order at the boundary).
\item If $P\in\Psi^m_\phi(X)$ is invertible as an operator $H^{s}_\phi(X)\to H^{s-m}_\phi(X)$ for some $s$ then its inverse is in $\Psi^{-m}_\phi(X)$.

\end{enumerate}

\begin{example}\label{ex:R^n}
 As a simple explicit example, consider $X=\overline{\R^n}$ with the Euclidean metric. 
We let $(r,y)\in (0,\infty)\times S^{n-1}$ be polar coordinates and set $x=\frac1r$ on $\R^n\setminus\{0\}$, so the Laplacian is $\Delta= D_r^2-i(n-1)r^{-1}D_r + r^{-2}\Delta_y = (x^2 D_x)^2 +i(n-1)x^3 D_x + x^2\Delta_y$, which has the form \eqref{eqn:phi operator}. 
Replacing $x^2 D_x$ by $\tau$ and $\Delta_y$ by $|\eta|^2$ one obtains
$\tau^2 + i(n-1)x\tau+|\eta|^2$, so the $\phi$-principal symbol and the normal family are both given by $\tau^2+|\eta|^2$. This is invertible for $(\tau,\eta)\neq0$, but not for $(\tau,\eta)=0$, so $\Delta$ is $\phi$-elliptic but not fully elliptic.

In contrast, the operator $\Delta+1$ has the same $\phi$-principal symbol but normal family $\tau^2+|\eta|^2+1$, hence is fully elliptic. Since it is invertible as an operator $H^2(\R^n)=H^2_\phi(X)\to L^2(\R^n)=H^0_\phi(X)$, items 5, 6 above imply that its inverse is in $\Psi_\phi^{-2}(X)$. This can also be seen explicitly: The Schwartz kernel of $(\Delta+1)^{-1}$ (with respect to Lebesgue measure, which is a positive $\phi$-density) is $K(p,p') = k(p-p')$ where $k=(|\xi|^2+1)\check{}\,$. So, for example, the rapid decay of $\beta_\phi^*K$ away from $\phif$ corresponds to the fact that $k$ is a Schwartz function, which implies rapid decay of $K(p,p')$ as $d(p,p')=|p-p'|\to\infty$
(see Remark \ref{rem:X2phi geometry}.2(a)).
Note that, in contrast, the Newton potential (corresponding to the inverse of $\Delta$) does not decay rapidly as $d(p,p')\to\infty$.
\end{example}
\begin{remark}[Connectedness II]
\label{rem:connected2}
If $\bdX$ is disconnected (cf.\ Remark \ref{rem:connected}) then the definition of $X^2_\phi$ needs to be modified slightly: Let $B_1,\dots,B_N$ be the connected components of $B$, and let $\calH_i = \phi^{-1}(B_i)$, so each $\calH_i$ is a union of components of $\bdX$. The b-blow-up \eqref{eqn:b-blowup} is replaced by
\begin{equation}
\label{eqn:b-blowup disconn}
\beta_b: X^2_b := [X^2;\bigcup_{i=1}^N (\calH_i)^2] \to X^2\,.
\end{equation}
Now $\bbf$ is a disjoint union of boundary hypersurfaces $\bbf_{GH}$, where $G, H$ vary over pairs of components of $\bdX$ contained in the same $\calH_i$.
The definition of the fibre diagonal $D_\phi$ and the $\phi$-blow-up \eqref{eqn:phi-blowup}, as well as the definition of $\Psi^m_\phi(X)$, remain the same, where now $\phif$ is a union of $\phif_{GH}$'s, one for each $\bbf_{GH}$. 
Note that $\diag_\phi$ meets only those $\phif_{GH}$ with $G=H$.%
\footnote{\label{fn:psido disconnected}This definition differs slightly from the usual one, e.g.\ as given in \cite[Sec.\ 3]{FriGriSch:CPFCO}, since it allows several components of $\bdX$ to be grouped together. Components in the same/different $\calH_i$'s should be considered close together/far apart: near (the lift of) $G\times H$ the kernel of an $A\in \Psi^m_\phi(X)$ is $O(x^\infty)$ if $G, H$ are in different $\calH_i's$, smooth to $\phif_{GH}$ if they are in the same $\calH_i$ and $G\neq H$, and smooth except for a conormal singularity at $\diag_\phi$ if $G=H$.
This greater flexibility is needed in Section \ref{sec:calderon}, see Footnote \ref{fn:connected Calderon}.
}
\end{remark}

\section{Resolvent and large time heat kernel}
\label{sec:resolvent heat}

The fundamental tool in our study of analytic torsion and index theory on 
$\phi$-manifolds is the microlocal analysis of the resolvent and the heat kernel. 
In this section we consider the resolvent of the Hodge Laplacian at low energy, and the heat kernel for large times. See Section \ref{subsec:heat short time} for the short time heat asymptotics.

Consider the open interior $\interior{X}$ of a $\phi$-manifold equipped with a $\phi$-metric $g_\phi$ and the de Rham complex $(\Omega^*(\interior{X}),d_*)$, where $d_k: \Omega^k(\interior{X}) \to \Omega^{k+1}(\interior{X})$ is the exterior differential acting on differential $k$-forms. Recall that the Hodge Laplace operator acting on differential $k$-forms is defined by 
$$
\Delta_k:= d_k^* d_k + d_{k-1} d_{k-1}^*,
$$
where the formal adjoints are defined with respect to $g_\phi$.
We omit the lower index $k$ whenever we consider the Laplacian acting on differential forms of all degrees. The Laplacians are positive operators, and we identify $\Delta$, $\Delta_k$ with their unique self-adjoint extensions in $L^2(\Omega^*(\interior{X}), g_\phi)$. \medskip

We can now proceed with the discussion of the resolvents $(\Delta_k+\kappa^2)^{-1}$ and $(\Delta+\kappa^2)^{-1}$ for $\kappa > 0$. 
We shall always abbreviate 
$$
\Lambda^k_\phi X:= \Lambda^k {}^\phi T^*X.
$$
The integral kernel of the resolvent $(\Delta_k+\kappa^2)^{-1},\kappa > 0$ is 
a section 
 of 
$(\Lambda^k_\phi X \otimes E)^* \boxtimes (\Lambda^k_\phi X \otimes E)$ pulled back to 
$(0,\infty)\times X\times X$ by the projection $(p,q,\kappa) \mapsto (p,q)$ onto the first two factors. 
We refer the reader to the careful discussion in \cite{gtv} concerning the precise choice of bundles and density factors involved, and explain here briefly only the blow-up of the base manifold 
$[0,\infty)\times X\times X$, necessary to turn the resolvent into a polyhomogeneous-conormal section. 
\medskip

Note that for $\kappa>0$ the operator $\Delta+\kappa^2$ is fully elliptic, so its inverse lies in $\Psi_\phi^{-2}(X)$, with smooth dependence on $\kappa$. So over $(0,\infty)\times X\times X$ our blown-up space will be simply $(0,\infty)\times X^2_\phi$. However, full ellipticity is lost at $\kappa=0$,\footnote{The normal family at $\tau=0$, $\eta=0$, i.e.\ $N(\Delta)(0;y,0)$, has kernel equal to the space of harmonic forms on the fibre $F_y$.} and this leads to singular behavior of the Schwartz kernel of $(\Delta+\kappa^2)^{-1}$ as $\kappa\to0$. Thus, the main point is to describe this behavior uniformly in terms of additional blow-ups over $\kappa=0$.

\medskip

This blow-up is described in detail in \cite[§ 6, § 7]{gtv} and is illustrated in
Figure \ref{fig:resolvent-blowup}, with the original space 
$[0,\infty)\times X\times X$ indicated with thick dotted (blue) 
coordinate axes in the background ($[0,\infty) \equiv \R_+$; only the $x$-coordinate of $X$ is shown, with $\partial X=\{x=0\}$).

\begin{figure}[h]
\centering
\begin{tikzpicture}[scale =0.55]

\draw[-] (-2,2.3)--(-2,5);
\draw[-] (-1,2)--(-1,5);
\draw[-] (1,2)--(1,5);
\draw[-] (2,2.3)--(2,5);
\draw[-] (-3,1)--(-6.5,-0.5);
\draw[-] (3,1)--(6.5,-0.5);
\draw[-] (-2.5,-0.5)--(-6.5,-2.3);
\draw[-] (2.5,-0.5)--(6.5,-2.3);

\draw (-1,2)--(-0.5,1.8);
\draw (1,2)--(0.5,1.8);
\draw (-0.5,1.8).. controls (0,2) and (0,2) .. (0.5,1.8);

\draw (-1,2).. controls (-1.1,1.9) and (-1.8,2.1) .. (-2,2.3);
\draw (1,2).. controls (1.1,1.9) and (1.8,2.1) .. (2,2.3);

\draw (-2.5,-0.5).. controls (-2,-1) and (-0.8,-1.2) .. (-0.5,-1.2);
\draw (2.5,-0.5).. controls (2,-1) and (0.8,-1.2) .. (0.5,-1.2);
\draw (-0.5,-1.2).. controls (-0.2,-1.5) and (0.2,-1.5) .. (0.5,-1.2);

\draw (-3,1).. controls (-2.2,0.7) and (-2.2,-0.2) .. (-2.5,-0.5);
\draw (3,1).. controls (2.2,0.7) and (2.2,-0.2) .. (2.5,-0.5);
\draw (-2,2.3).. controls (-2.4,2) and (-2.5,1.7) .. (-3,1);
\draw (2,2.3).. controls (2.4,2) and (2.5,1.7) .. (3,1);

\draw (-0.5,-1.2).. controls (0,-1) and (0,1.6) .. (-0.5,1.8);
\draw (0.5,-1.2).. controls (1,-1) and (1,1.6) .. (0.5,1.8);

\node at (-4,-0.3) {$\textup{lb}_{0}$};
\node at (4,-0.3) {$\textup{rb}_{0}$};
\node at (-4,2) {$\textup{lb}$};
\node at (4,2) {$\textup{rb}$};
\node at (0,-2) {$\textup{zf}$};
\node at (-1,1) {$\textup{bf}_{0}$};
\node at (0.41,0) {$\phi \textup{f}_0$};
\node at (0.5,3) {$\phi\textup{f}$};
\node at (-1.5,3) {$\textup{bf}$};

\draw[->, very thick, densely dotted, blue] (0.1,1)--(0.1,6);
\node[very thick, blue] at (1,6) {$\R_+$};
\draw[->, very thick, densely dotted, blue] (0.1,1)--(-7,-2);
\node[very thick, blue] at (-7.7,-1.5) {$X$};
\draw[->, very thick, densely dotted, blue] (0.1,1)--(7,-2);
\node[very thick, blue] at (7.7,-1.5) {$X$};

\end{tikzpicture}
 \caption{Resolvent blow-up space $X^{2}_{\kappa, \phi}$}
  \label{fig:resolvent-blowup}
\end{figure}

The blow-up is obtained as follows. First, we blow up 
the codimension $3$ corner $\{0\} \times \partial X \times \partial X$,
which creates a new boundary hypersurface $\textup{bf}_{0}$.
Then we blow up the (lifts of the) codimension $2$ corners
$\{0\}\times X \times \partial X$, $\{0\}\times \partial X \times X$
and $\R_+\times \partial X \times \partial X$, creating new 
boundary hypersurfaces $\textup{lb}_{0}$, $\textup{rb}_{0}$ and $\textup{bf}$, respectively. 
Next we blow up the (lift of the) fibre diagonal for all $\kappa$, i.e. $\R_+\times D_{\phi}$.
This defines a new boundary hypersurface $\phi\textup{f}$. Finally, 
the resolvent blow-up space $X^{2}_{\kappa, \phi}$ is obtained by one last blow-up
of the intersection of the (lifted) interior fibre diagonal, $\R_+\times \textup{diag}_{\phi,\mathrm{int}}$, (see \eqref{def diag phi int}) with $\textup{bf}_{0}$, which creates the
boundary hypersurface $\phi \textup{f}_0$. The composition of all these blow-down maps defines the total blow-down map
\[
\beta_{\kappa,\phi}\colon X_{\kappa,\phi}^2\to [0,\infty)\times X\times X.
\]
Our main theorem is that under additional conditions, the resolvent lifts to a polyhomogeneous section on the blown-up space
$X_{\kappa,\phi}^2$ with a conormal singularity along the lifted diagonal. 
More precisely, it lies in the so-called
split pseudo-differential calculus, where the word split refers to a splitting of the space of forms into the subspace of fibre-harmonic forms, where the singular behavior as $\kappa\to0$ occurs, and its orthogonal complement.
We omit the precise definition of this calculus, referring the reader to \cite[Definition 7.4]{gtv},
and instead just state somewhat informally the following theorem.

\begin{theorem}[{\cite[Theorem 7.11, Theorem 8.1]{gtv}}] \label{thm:IndexSets1}
Assume that $g_\phi$ is a perturbation of a product type $\phi$-metric $g_0$ on a $\phi$-manifold $X$, i.e.
\begin{equation}
g_{\phi} \restriction \calU = g_0 + w,\quad g_0=\frac{dx^{2}}{x^4} + \frac{\phi^{*}g_{B}}{x^{2}} + h,
\end{equation}
and assume
\begin{itemize}
\item the higher order term $w$ satisfies $\vert w \vert_{g_{0}} = O(x^3)$ as $x \to 0$,
\item $\phi: (\partial X, h + \phi^*g_B) 
\to (B, g_B)$ is a Riemannian submersion,
\item the base $B$ of the fibration $\phi: \partial X \to B$ has dimension $\dim B \geq 2$,
\item the additional spectral assumptions 1.4 and 1.5 in \cite{gtv} hold.
\end{itemize}
Then, for any integer $\sigma>0$, the lift of the Schwartz kernel of $(\Delta+\kappa^2)^{-\sigma}$ is a conormal distribution on
$X^{2}_{\kappa, \phi}$ vanishing to infinite order 
at rb, lb and bf, with the following index set bounds at some of the other faces 
\begin{align*}
\mathcal{E}_{\phi\textup{f}}\geq 0, \mathcal{E}_{\phi \textup{f}_0}\geq 
\min \, \{ \, 0, -2\sigma + (b+1) \, \}, \mathcal{E}_{\textup{bf}_0}\geq -2\sigma, 
\mathcal{E}_{\textup{zf}}\geq -2\sigma.
\end{align*}
\end{theorem}
See also \cite{KotRoc:LELRSFBO}. We continue under the assumptions of the theorem in our discussion of analytic torsion and the index theorem below. We use different assumptions in the discussion of the Calderón projector. \medskip

The idea behind the proof of the statement for $\sigma = 1$ is that (identifying operators and their Schwartz kernels) the lift
$\beta_{\kappa,\phi}^*(\Delta+\kappa^2)$ restricts to very explicit 
operators, called normal operators, at the various faces of $X^{2}_{\kappa, \phi}$, which then can be inverted to construct the resolvent iteratively.
The idea for $\sigma \geq 2$ is that the 
resolvent lies in  the calculus of operators with integral kernels being conormal distributions on $X^{2}_{\kappa, \phi}$ with the explicitly determined index set bounds. 
The construction uses the split b-$\phi$-pseudodifferential calculus constructed in  \cite{GriHun:PCLQLSS}, which simplifies and
generalizes the construction in \cite{vaillant}.
\medskip

In order to pass from the resolvent to the 
heat operator, recall that the resolvent and the heat operator are 
related by
\begin{equation}
e^{-t\Delta}=\frac{1}{2\pi i}\int_{\Gamma} e^{\xi t} (\Delta+\xi)^{-1}\, d\xi,
\label{eq:HeatOperator}
\end{equation}
with the usual contour $\Gamma = \Gamma(\theta,t)$ circling around the origin and  chosen to avoid the spectrum of $-\Delta$, see \cite[Figure\ 5]{JorgenBoris}.
The contour is chosen so that $\xi\in\Gamma(\theta,t)$ implies $\Re\xi<0$ or $\xi=t^{-1} e^{i\alpha}$ where $\alpha\not\in\pi+2\pi\Z$, so that $e^{t\xi}$ is bounded on $\Gamma$. This allows to relate the asymptotics of $e^{-t\Delta}$ for large $t$ 
to the asymptotics of  $(\Delta+\xi)^{-1}$ for $\xi=t^{-1} e^{i\alpha}=\kappa^2 e^{i\alpha}$ where $\kappa=t^{-1/2}$ is small. The latter is studied in the same way as for $\xi=\kappa^2$. The resulting
asymptotics of the heat kernel at the diagonal at large times are given by:%
\footnote{A minor technical point is that it is useful to avoid dealing with conormal singularities by considering the variant of \eqref{eq:HeatOperator}
\begin{equation*}
e^{-t\Delta}=\frac{(\nu-1)!}{(-t)^{(\nu-1)}}\frac{1}{2\pi i} \int_{\Gamma}  e^{\xi t} (\Delta+\xi)^{-\nu}\, d\xi,
\label{eq:ModifiedHeat}
\end{equation*}
If  
$\nu=\left[\frac{\dim X}{2}\right]+1
$
then the Schwartz kernel of $(\Delta+\xi)^{-\nu}$ is continuous.
}

\begin{theorem} \cite{JorgenBoris}
\label{thm:LongHeatKernelPoly}
The pointwise trace $\Tr(e^{-t\Delta})$ for large $t$
lifts to a polyhomogeneous function on the (lifted) diagonal of $X^{2}_{\kappa, \phi}$, where the variable $\kappa$ is replaced by $t^{-1/2}$,
with index set bounds 
$$\mathcal{E}_{\phi\textup{f}}\geq 0, \quad
\mathcal{E}_{\textup{zf}}\geq 0, \quad
\mathcal{E}_{\phi \textup{f}_0}\geq b+1,$$
at $\phi\textup{f}$, $\textup{zf}$ and $\phi \textup{f}_0$, respectively.
\end{theorem}

See Figure \ref{fig: cpt 0 heat space}. A full account of the asymptotics of the full heat kernel at all boundary faces of $X^{2}_{\kappa, \phi}$ can be obtained by following \cite{Sher}, who studied the same problem in the case of trivial fibres $F$.

\section{Analytic torsion}
\label{sec:torsion}

Analytic torsion is a central topological invariant, which is defined using spectral data on a compact Riemannian manifold. 
We start by recalling its definition. Let $(X,g)$ be a closed oriented Riemannian manifold of dimension $n$ equipped with a flat Hermitian
vector bundle $(E,\nabla,h)$. Consider the corresponding Hodge Laplacian $\Delta_k$
acting on $E$--valued differential forms $\Omega^k(X,E)$ of degree
$k$. Its unique self-adjoint extension is an operator with discrete spectrum $\sigma(\Delta_k)$
accumulating at $\infty$ according to Weyl's law. As a consequence, we can define its zeta-function as the following convergent sum
\begin{equation}
     \zeta(s,\Delta_k) :=  \sum_{\lambda\in \sigma(\Delta_k) \setminus \{0\}} m(\lambda) \, \lambda^{-s}, \ \Re( s) > \frac{n}{2},
     \label{eq:zeta}
\end{equation}
where $m(\lambda)$ is the geometric multiplicity of the eigenvalue $\lambda\in \sigma(\Delta_k)$.
The zeta-function is linked to the \emph{heat trace}
$\Tr \left( e^{-t\Delta_k}\right) $ via a Mellin transform
\begin{equation}
\zeta(s,\Delta_k) = \frac{1}{\Gamma(s)} \int_0^{\infty}
    t^{s-1} \left( \Tr \left( e^{-t\Delta_k}\right) - \dim \ker \Delta_k \right) \, dt.
    \label{eq:HeatTorsion}
\end{equation}
The short time asymptotic expansion of the
heat trace yields a meromorphic extension of $\zeta(s,\Delta_k)$, which is a priori a holomorphic function on $\Re( s) > \frac{n}{2}$,
to the whole complex plane $\C$ with at most simple poles and $s=0$ being a regular point.
This allows us to define the \emph{scalar analytic torsion} of the flat bundle $E$ by
\begin{equation}
T(X,E;g):= \exp\left(\, \frac{1}{2}\sum_{k=0}^{n}(-1)^k\cdot k \cdot 
\zeta'(0,\Delta_k)\right).
\label{eq:AnalyticTorsion}
\end{equation}

Below, it will become convenient to reinterpret the analytic torsion as a norm on the determinant line of the 
cohomology $H^*(X,E)$ with values in the flat vector bundle $E$. This is defined as 
\[\det H^*(X,E):= \bigotimes_{k=0}^m \left( \bigwedge^{b_k} H^k(X,E)\right)^{(-1)^k},\]
where $b_k=\dim H^k(X,E)$. The determinant line bundle has an induced $L^2$ structure $\| \cdot \|_{L^2}$ 
from the inclusion $H^k(X,E) \cong \ker \Delta_{k} \subset L^2\Omega^k(X,E)$, where the $L^2$-inner product
is induced by the Riemannian metric $g$ and the Hermitian metric $h$. 
The \emph{analytic torsion norm}, also called the Ray-Singer metric or the Quillen metric, is then defined as
\begin{equation}
\| \cdot \|_{(X,E,g)}^{RS}:= T(X,E;g)\| \cdot \|_{L^2}.
\label{eq:AnalyticTorsionNorm}
\end{equation}
By an argument of Ray and Singer, this is actually independent of $g$. \medskip

In order to extend this construction to a non-compact or singular manifold, with compactification $X$, we replace $H^*(X,E)$ by the reduced $L^2$-cohomology $H^*_{(2)}(X,E)$, which is  isomorphic to the space of $L^2$-harmonic forms. The time integral in \eqref{eq:HeatTorsion} is replaced by a renormalized integral and the heat trace ($e^{-t\Delta_k}$ is no longer trace class) is replaced by a renormalized heat trace. Its definition and asymptotics is the central analytic difficulty that requires Theorem \ref{thm:LongHeatKernelPoly}.  Our first main result in \cite{JorgenBoris} is precisely the statement that it is indeed possible.

\begin{theorem}\cite[Theorem 1.3]{JorgenBoris} Let $(X,g_\phi)$ be a $\phi$-manifold and $(E,h)$ a flat vector bundle over $X$.
We impose the assumptions of Theorem \ref{thm:IndexSets1}. Then the following statements hold: 
\begin{enumerate}
\item the renormalized analytic torsion $T(X,E;g_\phi)$ is well-defined; 
\item if $n=\dim(X)$ is odd, then the torsion norm $\lvert|\cdot |\rvert_{(X,E,g_\phi)}^{RS}$ on the determinant line 
$\det H^*_{2}(X;E)$ of the reduced $L^2$-cohomology is invariant under perturbations of $g_\phi$ 
of the form $\delta g_\phi= h$ with 
$\int_X \vert h\vert_{g_\phi} \, d\text{Vol}_{\phi}<\infty$. 
\end{enumerate}
\end{theorem}

A second interesting observation is the fact that the renormalized analytic torsion on a $\phi$-manifold still admits a cut and paste property. An interesting application is then that $T(X,E;g_\phi)$ appears as a correction term in a singular degeneration process. Let us explain the result informally, referring the reader to \cite{JorgenBoris} for details. We consider a sequence $(K_\varepsilon, g_\varepsilon)$ of compact closed manifolds that degenerate to a space $(\Omega, g_\omega)$ with a possibly non-isolated conical (wedge) singularity. This singular degeneration is constructed by gluing rescaled parts of a $\phi$-manifold $(\interior{X},g_\phi)$ to a fixed compact manifold with boundary, where we assume that the boundary fibration is trivial, namely $\partial X = B \times F$. We obtain the following result. 

\begin{theorem}\cite[Theorem 1.6]{JorgenBoris} Consider the singular degeneration as described above. Impose the assumptions of Theorem \ref{thm:IndexSets1}, assume $\dim F$ is even. and also one of the following two assumptions:
\begin{enumerate}
\item either $E$ is acyclic:  $H^k(\partial X,E)=0$ for all degrees $0\leq k \leq \dim(\partial X)$,
\item or $H^k(\partial X,E)=0$ for $1\leq k \leq \dim(\partial X)-1 $ and 
$E$ is a trivial vector bundle over $X$.
\end{enumerate}
Then, for any $\varepsilon>0$ there is a canonical isomorphism of determinant lines 
\[
\Phi: \det H^*_{(2)} (X,E) \otimes \det H^*_{(2)} (\Omega, E) \to \det H^* (K_{\varepsilon},E) ,\]
which is an isometry of the (renormalized) analytic torsion norms. Moreover
\begin{equation}\label{gluing-corr2}
\begin{split}
\log T(K_\varepsilon, E; g_\varepsilon) &=
\log T(X,E;g_\phi) + \log T(\Omega,E;g_\omega) \\ &+ 
\log\frac{\lvert |{\alpha}|\rvert_{L^2(X,E;g_\phi)} \cdot \lvert |{\beta}|\rvert_{L^2(\Omega,E;
g_\omega)}}{\lvert |{\Phi (\alpha \otimes \beta )|\rvert}_{L^2(K_\varepsilon, E; g_\varepsilon)}}.
\end{split}
\end{equation}
\end{theorem}

\section{Index theory and short time heat kernel}
\label{sec:index}

We can refine the description of the resolvent and heat kernel of the Laplacian described in Section \ref{sec:resolvent heat} to try to obtain an index theorem for this class of manifolds. In the philosophy of Melrose \cite[p. 1]{mel:aps}, such a result corresponds to the Atiyah-Singer index theorem in \textit{the category of $c$-$\phi$-manifolds}.

\subsection{The Atiyah-Singer index theorem and the heat kernel}

Let us formulate briefly what we mean by index theorem. For a historical overview of the subject, a good reference is \cite{freed}.
The classical setup is the following: consider vector bundles $E,E'\to X$ over a compact manifold without boundary. Let $P:C^\infty(X;E)\to C^\infty(X;E')$ be a differential operator of order $m$ between smooth sections of these bundles. If $P$ is also elliptic (meaning its principal symbol is invertible over $T^*X \setminus \{0\}$), compactness ensures that its extension as an operator between Sobolev spaces $H^s(E)\to H^{s-m}(E')$ is unique and Fredholm for all $s\in\R$. The index theorem of Atiyah and Singer \cite{atisin:1963} computes its Fredholm index in terms of characteristic classes:
$$\text{ind} \left(P\right) = (-1)^n \int_{TX}\text{Ch}(P)\text{Todd}(TX).$$
Notice how the left hand side is analytical in nature, while the right hand side is topological. This result is a generalization of several other important statements and permeates the areas of Analysis, Geometry and Topology since the second half of the 20th century. It has been looked at through several lenses, like cobordism and K-theory. We will focus on one of them as a promising strategy to generalize the result to our context: the so-called \textbf{heat kernel method}.

This is based on the fact that in this case of closed manifolds\footnote{Whenever the vector bundles have a hermitian structure (and the base manifold $X$ a Riemannian metric) that allows us to define the adjoint $P^*:C^\infty(X;E')\to C^\infty(X;E)$.} \footnote{In the case of Dirac operators $P=\dirac$ on even dimensional manifolds, the supertrace corresponds to the trace graded with respect to the superbundle structure on the Clifford module \cite[p. 120]{bgv}, hence its name.}:
$$\text{ind} \left(P\right) = \lim_{t\to\infty} \text{Str } e^{-t\hat{P}^2}, \quad \text{where} \quad \text{Str } e^{-t\hat{P}^2} = \text{Tr } e^{-tP^*P} - \text{Tr } e^{-tPP^*},$$
$$\hat{P} = \begin{pmatrix}
    0 &  P^* \\
    P & 0
    \end{pmatrix}: C^\infty(X;E)\oplus C^\infty(X;E') \longrightarrow C^\infty(X;E)\oplus C^\infty(X;E'),$$
so one can (hope to) obtain an index formula in more general geometries by understanding the asymptotics of the supertrace of the heat operator. It is useful to write:
$$\lim_{t\to\infty} \text{Str } e^{-t\hat{P}^2} = \lim_{t\to 0} \text{Str } e^{-t\hat{P}^2} + \int_0^\infty \partial_t \text{Str } e^{-t\hat{P}^2} dt.$$
We will call this the McKean-Singer formula. By Lidskii's theorem:
$$\text{Tr } A = \int_X K_A(x,x),$$
so this amounts to understanding the asymptotics of the heat kernel $h(t,x,x')$ of $\hat{P}^2$ at the diagonal $\R_+\times\text{diag}_X\subset \R_+\times X^2$.

In \cite{aps:1975} the case of a compact manifold with boundary and $P=\dirac^+$ a Dirac operator is examined. The metric is assumed to be of product type in a neighborhood of the boundary. To set up a Fredholm problem for the operator, a key step is recognizing the existence of a topological obstruction to the existence of local boundary conditions giving Fredholmness, and accordingly finding appropiate global (now called APS) boundary conditions. This corresponds to attaching an infinite half-cylinder of cross section $\partial X$ at the boundary and restricting the domain to $L^2$-sections along this \textit{cylindrical end}. They obtain the following index formula\footnote{Instead of $\hat{P}$ and the bundles $E$ and $E'$, we now have a decomposition coming from the $\mathbb{Z}_2$-grading of the Clifford bundle $E\to X$ associated to the Dirac operator:
$$\dirac = \begin{pmatrix}
    0 &  \dirac^- \\
    \dirac^+ & 0
    \end{pmatrix}: C^\infty(X;E^+)\oplus C^\infty(X;E^-) \longrightarrow C^\infty(X;E^+)\oplus C^\infty(X;E^-).$$
    $\dirac^+_{\text{APS}}$ refers to the Dirac operator with APS boundary conditions.}:
$$\text{ind} \left(\dirac^+_{\text{APS}}\right) + \frac{1}{2}\dim \ker \dirac_{\text{APS}, \partial X} = \int_X AS(x)dx - \frac{1}{2}\eta(\dirac_{\text{APS},\partial X}),$$
where $AS$ is effectively the same integrand as in the closed case\footnote{$AS$ is a priori a cohomology class, but the analysis of the heat kernel asymptotics provides us with a differential form representative, making sense of the integral.}, coming from the short time asymptotics of the heat kernel on the doubled manifold, and $\eta(\dirac_{\text{APS},\partial X})$ ``counts'' the number of positive eigenvalues minus the number of negative ones for the induced Dirac operator on the boundary $\dirac_{\text{APS},\partial X}$.

On the other hand, \cite{mel:aps} looks directly at manifolds with cylindrical ends and computes the different terms in the McKean-Singer formula, even in the non-Fredholm case (which corresponds to non-invertible induced boundary operator). For that, Melrose considers Dirac operators and lets them act on weighted Sobolev spaces, for which Fredholmness only fails if the weight belongs to the discrete set $-\text{Spec}(\dirac_{\partial X})$. The (generalized) index for a non-Fredholm weight $s$ is:
$$\text{ind}_s \left(\dirac^+\right) \coloneqq\lim_{\varepsilon\to 0}\frac{\text{ind}_{s - \varepsilon}\left(\dirac^+\right) + \text{ind}_{s + \varepsilon} \left(\dirac^+\right)}{2}.$$ 
In the case of a non-invertible boundary operator, this yields a way of computing the original index (i.e for weight $s = 0$), which reproduces the Atiyah-Patodi-Singer result.

With that, Melrose not only solves the case of cylindrical ends, but also establishes a general philosophy to address index theory of Dirac operators in singular spaces, which can be in particular applied to $c$-$\phi$-metrics.

A few comments about the McKean-Singer formula are here in place:
\begin{itemize}
    \item As the treatment in \cite[p. 10-12]{mel:aps} shows, even if $P$ is not Fredholm, a priori we can still compute the long time limit of the supertrace, so Fredholmness is not required for the method to work, permitting more general results than other approaches.
    \item It can however be very hard to compute the specific contributions to the formula, even in the closed case (cf. \cite[p. 29]{freed})
    . We will therefore have to restrict ourselves to the case of $P$ being a Dirac-type operator, as in \cite{mel:aps}, where adapting an argument of Getzler \cite{getzler} will allow us to recover the correct terms for the McKean-Singer formula. For a Clifford bundle $E \to X$ over a closed spin manifold with associated Dirac operator $\dirac_E$, the index formula one obtains via the heat kernel method is \cite[Theorem 4.3]{bgv}:
    $$\text{ind} \left(\dirac^+_E\right) = \int_X \hat{\text{A}} (TX) \text{Ch}(E),$$
    where the characteristic classes in the right hand side are called \textit{A-hat genus} and \textit{Chern character} \cite[p. 47-48, 142-143]{bgv}. In the Atiyah-Patodi-Singer case (with boundary), the formula looks like \cite[p. 10-12]{mel:aps}:
    \begin{align*}
        & \text{ind} \left(\dirac^+_E\right) = \int_X \hat{\text{A}} (TX) \text{Ch}(E) -\underbrace{\frac{1}{2\sqrt{\pi}}\int_0^\infty \text{Tr}\left(\dirac_{E,\partial X} e^{-t\left(\dirac_{E,\partial X}\right)^2}\right)\frac{dt}{t^{1/2}}}_{\frac{1}{2}\eta(\dirac_{E,\partial X})}, \\
    & \text{where} \quad \text{ind} \left(\dirac^+_E\right) = \dim\ker \left(\dirac^+_E\right) - \dim\ker_- \left(\dirac^-_E\right) + \frac{1}{2}\dim\ker \dirac_{E,\partial X} \\
    & \text{and} \quad \ker_- \left(\dirac^-_E\right) = \bigcap_{\varepsilon > 0} \ker \restr{\left(\dirac^-_E\right)}{x^{-\varepsilon}L^2}.
    \end{align*}
    
    \item There is in general no reason for the Lidskii integral of the kernel at the diagonal, which appears in each member of the formula through the supertrace, to converge. We can solve this by rather considering their \textit{renormalized} version \cite[p. 22-23]{alb:renorm}:
    $$\lim_{t\to\infty} \leftindex^R {\text{Str }} e^{-t\hat{P}^2} = \lim_{t\to 0} \leftindex^R {\text{Str }} e^{-t\hat{P}^2} + \int_0^\infty \partial_t \leftindex^R {\text{Str }} e^{-t\hat{P}^2} dt,$$
    where
$$\leftindex^R {\text{Tr }} A = \leftindex^R {\int}_X K_A(p,p) \quad \text{and} \quad \leftindex^R {\int}_X f = \underset{z = 0}{\text{FP}} \int_X x^zf $$
    when $f$ is polyhomogeneous conormal\footnote{The integral $h(z)=\int_X x^z f$ is a priori defined for $\textup{Re}(z)$ large, but has a meromorphic continuation to $z\in\mathbb{C}$. $\underset{z=0}{\textup{FP}} \; h(z)$ denotes the $z^0$-term in its Laurent expansion at $z=0$.}, for a choice of boundary defining function $x$. Here, $\leftindex^R{\int_X f}$ is called \textit{the renormalized integral of $f$}. Alternatively, one could integrate $f$ over the truncated manifold $\{x\geq\varepsilon\}$ and compute the coefficient of $\varepsilon^0$, obtaining the same result in the cases of interest to us.

    This \textit{generalized} McKean-Singer formula can be applied to any context, as long as the heat kernel is shown to be polyhomogeneous. As in the formulas above, the large time limit will correspond to the index of the operator, when it is Fredholm, or some generalized version of it otherwise. The Atiyah-Singer integrand appears in the short time limit, but this limit can also add more contributions, which together with the integral term in $t\in(0,\infty)$ are usually boundary terms like the $\eta$-invariant. 
    
    As a consequence, the first step in the heat kernel method will be to blow up the space $\R_+\times X^2$, where the heat kernel lives, to resolve its singularities i.e. make it (lift to) a polyhomogeneous distribution. This encodes the complexity of the asymptotics of the heat kernel geometrically. One may call this step \textit{finding the heat blow-up space}.

\end{itemize}

\subsection{The heat blow-up spaces}
\label{subsec:heat short time}

The heat blow-up space for short times is constructed, and the heat kernel is shown to be polyhomogeneous on it, for $c=0$ in \cite[\S 5, Corollary 7.2]{tave}, for $c=1$ in \cite[\S 4.1, Theorem 4.11]{vaillant} and for $c=2$ in \cite[\S 3.2.1]{ash:thesis}. \footnote{Brownian motion was used in \cite{simek:heatcusp} to study the (first three terms of the) short time Dirichlet heat trace asymptotics for planar domains and surfaces of revolution with cusps of any sharpness $\gamma>1$ (and provide bounds on the Neumann heat trace in the planar case); cf. Footnote \ref{fn:other cusps} and Figure \ref{fig:examples} left. In our terminology the $\gamma=2$ case of this corresponds to $c=2$, trivial base, and fibre equal to a ball (in particular the fibre has a boundary).
}

The guiding principle in all of these constructions is to make the function $\frac{d(p,p')^2}{t}$ lift to a polyhomogeneous function with as few blow-ups as possible, see \cite[\S 3.2.1]{ash:thesis}. Here, $d$ denotes the Riemannian distance. Compare Remark \ref{rem:X2phi geometry}.

\begin{figure}[!ht]
\centering
\resizebox{.5\textwidth}{!}{%
\begin{circuitikz}
\tikzstyle{every node}=[font=\LARGE]

\draw (20,14.5) to[short] (20,23.75);
\draw (11.25,14) to[short] (11.25,23);
\draw [short] (11.25,14) .. controls (11.75,13.5) and (12.5,12.75) .. (13.75,12.75);
\draw [short] (19.25,13.75) .. controls (20,14) and (20,14.25) .. (20,14.5);
\draw [->, >=Stealth] (16.25,20.5) -- (16.25,25.25);
\draw [->, >=Stealth] (20,14.5) -- (23.75,13.25);
\draw [->, >=Stealth] (11.25,14) -- (8.25,11.25);
\draw [short] (15.75,12.25) .. controls (16.5,13.5) and (17.5,13.5) .. (18.25,12.75);
\draw [short] (13.75,12.75) .. controls (14.25,12.5) and (15,12.25) .. (15.75,12.25);
\draw [short] (18.25,12.75) .. controls (18.75,13) and (19,13.5) .. (19.25,13.75);
\draw [short] (15.75,12.25) -- (18.75,7.75);
\node [font=\Huge] at (16,15.75) {};
\node [font=\Huge] at (18.25,10.75) {\fontsize{30}{28}\selectfont tf};
\node [font=\Huge] at (17,25.5) {\fontsize{30}{28}\selectfont $t$};
\node [font=\Huge] at (24.5,13.5) {\fontsize{30}{28}\selectfont $x'$};
\node [font=\Huge] at (7.8,11.3) {\fontsize{30}{28}\selectfont $x$};
\draw [short] (18.75,7.75) .. controls (19.5,9) and (20.5,9) .. (21.25,8.25);
\draw [short] (18.25,12.75) -- (21.25,8.25);
\node [font=\Huge] at (16.5,17.25) {\fontsize{30}{28}\selectfont $\phi$f};
\node [font=\Huge] at (12.5,18) {\fontsize{30}{28}\selectfont bf};
\draw (13.75,12.75) to[short] (13.75,22.5);
\draw (19.25,13.75) to[short] (19.25,23.25);
\end{circuitikz}
}%
\caption{Heat blow-up space for $0$-$\phi$-metrics (compare \cite[Figure 4]{tave}). The face $\phi\textup{f}$ here corresponds to $\phi\textup{f}$ in Figure \ref{fig:resolvent-blowup}.}
\label{fig: 0 heat space}
\end{figure}

\begin{figure}[!ht]
\centering
\resizebox{.5\textwidth}{!}{%
\begin{circuitikz}
\tikzstyle{every node}=[font=\LARGE]

\draw (20,14.5) to[short] (20,23.75);
\draw (11.25,14) to[short] (11.25,23);
\draw [short] (11.25,14) .. controls (11.75,13.5) and (12.5,12.75) .. (13.75,12.75);
\draw [short] (19.25,13.75) .. controls (19.75,14) and (20,14.25) .. (20,14.5);
\draw [->, >=Stealth] (16.25,20.5) -- (16.25,25.25);
\draw [->, >=Stealth] (20,14.5) -- (23.75,13.25);
\draw [->, >=Stealth] (11.25,14) -- (8.25,11.25);
\draw [short] (15.75,12.25) .. controls (16.5,13.5) and (17.5,13.5) .. (18.25,12.75);
\draw [short] (13.75,12.75) .. controls (14.25,12.5) and (15,12.25) .. (15.75,12.25);
\draw [short] (18.25,12.75) .. controls (18.75,13) and (19,13.5) .. (19.25,13.75);
\draw [short] (15.75,12.25) -- (18.75,7.75);
\node [font=\LARGE] at (16,15.75) {};
\node [font=\LARGE] at (18.25,10.75) {\fontsize{30}{28}\selectfont tf};
\node [font=\LARGE] at (17,25.5) {\fontsize{30}{28}\selectfont $t$};
\node [font=\LARGE] at (24.5,13) {\fontsize{30}{28}\selectfont $x'$};
\node [font=\LARGE] at (7.8,11.3) {\fontsize{30}{28}\selectfont $x$};
\draw [short] (18.75,7.75) .. controls (19.5,9) and (20.5,9) .. (21.25,8.25);
\draw [short] (18.25,12.75) -- (21.25,8.25);
\draw [short] (13.75,12.75) .. controls (14.25,16.5) and (18.25,17.25) .. (19.25,13.75);
\node [font=\LARGE] at (16.25,18.75) {\fontsize{30}{28}\selectfont bf};
\node [font=\LARGE] at (16.5,14.25) {\fontsize{30}{28}\selectfont $\phi$tf};
\end{circuitikz}
}%
\caption{Heat blow-up space for $1$-$\phi$-metrics (compare \cite[Figure 5]{vaillant}).}

\label{fig: 1 heat space}
\end{figure}

\begin{figure}[!ht]
\centering
\resizebox{.5\textwidth}{!}{%
\begin{circuitikz}
\tikzstyle{every node}=[font=\LARGE]

\draw [short] (11.25,14) .. controls (11.75,13.5) and (12.5,12.75) .. (13.75,12.75);
\draw [short] (19.25,13.75) .. controls (20,14) and (20,14.25) .. (20,14.5);
\draw [->, >=Stealth] (15,18.25) -- (15,24.75);
\draw [->, >=Stealth] (20,14.5) -- (23.75,13.25);
\draw [->, >=Stealth] (11.25,14) -- (8.25,11.25);
\draw [short] (15.75,12.25) .. controls (16.5,13.5) and (17.5,13.5) .. (18.25,12.75);
\draw [short] (13.75,12.75) .. controls (14.25,12.5) and (15,12.25) .. (15.75,12.25);
\draw [short] (18.25,12.75) .. controls (18.75,13) and (19,13.5) .. (19.25,13.75);
\draw [short] (15.75,12.25) -- (18.75,7.75);
\node [font=\LARGE] at (16,15.75) {};
\node [font=\LARGE] at (18.25,10.75) {\fontsize{30}{28}\selectfont tf};
\node [font=\LARGE] at (16,25) {\fontsize{30}{28}\selectfont $t$};
\node [font=\LARGE] at (24.5,13) {\fontsize{30}{28}\selectfont $x'$};
\node [font=\LARGE] at (7.8,11.3) {\fontsize{30}{28}\selectfont $x$};
\draw [short] (18.75,7.75) .. controls (19.5,9) and (20.5,9) .. (21.25,8.25);
\draw [short] (18.25,12.75) -- (21.25,8.25);
\draw [short] (20,14.5) .. controls (19,18.25) and (14.25,21.25) .. (11.25,14);
\draw [short] (19.25,13.75) .. controls (18,15.5) and (15.5,16.75) .. (13.75,12.75);
\node [font=\LARGE] at (15.5,16.5) {\fontsize{30}{28}\selectfont btf};
\node [font=\LARGE] at (16.5,14) {\fontsize{30}{28}\selectfont $\phi$tf};
\end{circuitikz}
}%
\caption{Heat blow-up space for $2$-$\phi$-metrics (compare \cite[Example 3.25, $(k,q)=(0,2)$]{ash:thesis}).}

\label{fig: 2 heat space}
\end{figure}

From the heat blow-up spaces, one can read off which faces contribute to each of the terms in the McKean-Singer formula (for each case $c\in\{0,1,2\}$ only the faces actually present in the respective heat space should be taken):
\begin{itemize}
    \item The faces tf, $\phi$tf and btf will contribute to the short time limit ($t\to 0$).
    \item $\phi$f and bf correspond to the integral term (heuristically, no contribution is expected from this term in the case $c=2$ due to the lack of a blown-up face at $t>0$).
    \item The blow-up spaces so constructed only deal with the heat kernel at finite times. One could extend them through compactification at large times by introducing a local coordinate $\kappa = \frac{1}{\sqrt{t}}$ and resolving the heat kernel for $\kappa \to 0$, see Figure \ref{fig: cpt 0 heat space}. To do this, as mentioned in Section \ref{sec:resolvent heat}, one usually studies the low energy resolvent of the Laplacian \cite[\S 6.6-6.9, \S 7.7-7.8]{mel:aps} \cite{gtv} and applies the functional calculus formula (\ref{eq:HeatOperator}). In this way, the heat blow-up space at large times corresponds to the blow-up space for the resolvent at low energy (Figure \ref{fig:resolvent-blowup}), and the faces at $\kappa=0$ will contribute to the asymptotics as $t\to\infty$. Another approach is to use the resolvent of the Dirac operator and apply Stone's formula \cite[(62)]{vaillant}:
    $$e ^{-t\dirac^2} = \frac{1}{2\pi i}\int_\R e^{-tr^2} \lim_{\varepsilon\to 0}\left(\left(\dirac - r - i\varepsilon\right)^{-1} - \left(\dirac - r + i\varepsilon\right)^{-1}\right)dr.$$

\end{itemize}

\begin{figure}[!ht]
\centering
\resizebox{.5\textwidth}{!}{%
\begin{circuitikz}
\tikzstyle{every node}=[font=\LARGE]

\draw [short] (11.25,14) .. controls (11.75,13.5) and (12.5,12.75) .. (13.75,12.75);
\draw [short] (19.25,13.75) .. controls (20,14) and (20,14.25) .. (20,14.5);
\draw [->, >=Stealth] (20,14.5) -- (23.75,13.25);
\draw [->, >=Stealth] (11.25,14) -- (8.25,11.25);
\draw [short] (15.75,12.25) .. controls (16.5,13.5) and (17.5,13.5) .. (18.25,12.75);
\draw [short] (13.75,12.75) .. controls (14.25,12.5) and (15,12.25) .. (15.75,12.25);
\draw [short] (18.25,12.75) .. controls (18.75,13) and (19,13.5) .. (19.25,13.75);
\draw [short] (15.75,12.25) -- (18.75,7.75);
\node [font=\LARGE] at (16,15.75) {};
\draw [short] (18.75,7.75) .. controls (19.5,9) and (20.5,9) .. (21.25,8.25);
\draw [short] (18.25,12.75) -- (21.25,8.25);
\draw [short] (11.25,14) -- (11.25,25.25);
\draw [short] (13.75,12.75) -- (13.75,24.75);
\draw [short] (19.25,13.75) -- (19.25,24.75);
\draw [short] (20,14.5) -- (20,25.75);
\draw [short] (20,25.75) .. controls (19.75,25.75) and (19.5,25.5) .. (19.25,24.75);
\draw [short] (13.75,24.75) .. controls (13,25.5) and (12.25,25.75) .. (11.25,25.25);
\draw [short] (11.25,25.25) .. controls (11,25.5) and (10.75,26) .. (10.75,26.25);
\draw [short] (10.75,26.25) .. controls (11.5,26.75) and (12,27.5) .. (11.5,28.5);
\draw [short] (13.75,24.75) .. controls (14.25,25.25) and (14.75,25.25) .. (15.5,25.25);
\draw [short] (15.5,25.25) .. controls (16,24.5) and (17.25,24.25) .. (18,25.25);
\draw [short] (18,25.25) .. controls (18.25,25.25) and (18.75,25.25) .. (19.25,24.75);
\draw [short] (20,25.75) .. controls (20.25,26.25) and (20.25,26.5) .. (20.25,26.75);
\draw [short] (20.25,26.75) .. controls (19.5,27.25) and (19,27.75) .. (19,28.75);
\draw [short] (14.25,28.75) .. controls (14,29.25) and (16.25,29.5) .. (16.75,28.75);
\draw [short] (15.5,25.25) .. controls (15.5,27.25) and (15.5,27.75) .. (14.25,28.75);
\draw [short] (18,25.25) .. controls (18,27.75) and (17.25,28.25) .. (16.75,28.75);
\node [font=\LARGE] at (8,27.25) {};
\draw [short] (10.75,26.25) -- (6.5,28.25);
\draw [short] (11.5,28.5) -- (7.25,30.5);
\draw [short] (20.25,26.75) -- (24.5,29.25);
\draw [short] (19,28.75) -- (23.25,31.25);
\draw [short] (11.5,28.5) .. controls (12,29) and (13.25,29) .. (14.25,28.75);
\draw [short] (19,28.75) .. controls (18.5,29) and (17,29) .. (16.75,28.75);
\node [font=\LARGE] at (18.5,10.5) {\fontsize{30}{28}\selectfont tf};
\node [font=\LARGE] at (16.5,18) {\fontsize{30}{28}\selectfont $\phi$f};
\node [font=\LARGE] at (16.5,27) {\fontsize{30}{28}\selectfont $ \phi\text{f}_0$};
\node [font=\LARGE] at (16,31.5) {\fontsize{30}{28}\selectfont zf};
\end{circuitikz}
}%
\caption{Compactified heat blow-up space for $0$-$\phi$-metrics (compare Figure \ref{fig: 0 heat space} and \cite[Figure 9]{gtv}). See also \cite{sher:ac_heat} for the case with trivial fibres.}

\label{fig: cpt 0 heat space}
\end{figure}

\subsection{When is $\dirac$ self-adjoint and Fredholm?}

A word on self-adjointness and Fredholmness of the respective Dirac operators:
\begin{itemize}
    \item For incomplete spaces (as it is the case for $c=2$), (essential) self-adjointness is not a given and there could be many self-adjoint extensions. The minimal and maximal closed extensions are
    \begin{align*}
        \mathcal{D}_{\text{min}} & = \{u\in L^2(X, E) : \exists(u_n)\subset C_0^\infty(\mathring{X}) \quad \text{so that} \quad u_n \to u \quad \text{and} \quad \dirac u_n \to \dirac u \}, \\
        \mathcal{D}_{\text{max}} & = \{u\in L^2(X, E) : \dirac u \in L^2(X,E)\}. 
    \end{align*}
    The quotient $\mathcal{D}_{\text{max}}/\mathcal{D}_{\text{min}}$ could even be infinite dimensional \cite[p. 5]{albgr:incedge}. To obtain uniqueness of the self-adjoint extension one usually imposes a so-called \textit{geometric Witt condition} on the vertical family $(\dirac_{F_y})_{y\in B}$ of Dirac operators induced on the boundary, namely:
    \begin{equation}\label{eqn:geometric Witt}
        \text{Spec}\left(\dirac_F\right) \cap \left(-\frac{1}{2}, \frac{1}{2}\right) = \emptyset,    
    \end{equation}
    in the settings where the (fibres at the) boundary undergo conical degeneration \cite[(3.2)]{chou} \cite[\S 3.1]{albgr:incedge} \cite[\S 2.2]{albgr:pseudo}. In some cases where the degeneration is cusp-like (of any sharpness, cf. Footnote \ref{fn:other cusps}), invertibility is enough \cite[Theorem 1]{grswo} \cite[Theorem 1.1]{liu:thesis} \cite[Theorem 5.3]{lespeye}). These conditions ensure that $\mathcal{D}_{\text{min}} = \mathcal{D}_{\text{max}}$ and have to do with (but the geometric Witt condition is stronger than) the vanishing of middle degree intersection homology (so-called \textit{Witt condition}).

    With completeness (for us, if $c\in\{0, 1\}$), (essential) self-adjointness comes for free \cite{gaf}. \footnote{We need the self-adjointness of the Dirac operator to carry out our heat kernel analysis (to have real eigenvalues, functional calculus, etc.) and consequently $\text{ind}(\dirac) = 0$. This is not to be confused with the non-trivial $\text{ind}(\dirac^+)$ that the index formula actually computes.}
    \item Fredholmness in suitable Sobolev spaces for $c = 0$ corresponds to the invertibility of the vertical boundary family, i.e. having a fully elliptic $\dirac$ (compare \cite[Theorem 1]{MM98}; see however \cite{GriHun:PCLQLSS}). For $c = 1$, it is the invertibility of a horizontal operator that controls the Fredholmness \cite[p. 2]{vaillant}. For $c= 2$, one could expect it to come for free once we have a self-adjoint operator. 
    
    When the operator is Fredholm, the large time limit of the supertrace is just the usual $L^2$-index (read off from the interior face at infinite time zf).
    
    Nevertheless, one could also work without Fredholmness and obtain a \textit{generalized index}, like in \cite[p. 11]{mel:aps} \cite[p. 63]{vaillant}. When the boundary is the total space of a fibration, one can handle the arising complexity by assuming that the kernel of the boundary family forms a vector bundle $\mathcal{K}\to B$ over the base \cite[(39)]{vaillant}, in a similar fashion as the treatment of the families index theorem (cf. \cite[p. 313]{bgv}). Then the horizontal operator mentioned above for $c = 1$ is precisely $\dirac_{\mathcal{K}, B}$.
\end{itemize}
    
\subsection{Rescaling and contributions in the $c$-$\phi$ index theorems}

Once the blow-up spaces are clear, we ``just'' need to figure out what terms in the asymptotic expansion at each face contribute to the McKean-Singer formula and compute them, with the help of Getzler rescaling. What we mean by that is roughly the following: for each $t\geq 0$ the heat kernel is a section of the bundle\footnote{Not to confuse with the bundle $\text{hom}(E)\to X$ with fibres $\text{hom}(E)_x=\text{hom}(E_x,E_x)$.} $\text{Hom}(E) \to X^2$, whose fibres are $\text{Hom}(E)_{(x,x')} = (E^* \boxtimes E)_{(x,x')}=\text{hom}(E_x, E_{x'})$. At the diagonal, where the interesting (for us) part of the kernel lives (Lidskii), $\restr{\text{Hom}(E)}{\text{diag}_X}\cong \restr{\text{hom}(E)}{X}\cong \mathbb{C}l(\leftindex^{c-\phi}TX) \otimes \text{hom}_{\mathbb{C}l}(E)$. \footnote{Here, $\leftindex^{c-\phi} TX = x^{-c}\left(\leftindex^\phi TX\right)$.} Now, the Clifford algebra has a natural filtration (by Clifford degree) inherited from $\Lambda \leftindex^{c-\phi}{T^*X}$, which gives us a filtration for $\restr{\text{Hom}(E)}{\text{diag}_X}$ by the previous isomorphism. This filtration can be first pulled back along the projection $p:\mathbb{R}_+\times X^2\to X^2$ and then to the heat blow-up space along $\beta:HX^2_\phi \to \mathbb{R}_+\times X^2$. The faces intersecting the lifted diagonal are contained in $\beta^{-1}\left(\mathbb{R}_+\times \text{diag}_X\right)$, so we obtain a filtration of $\left(\beta\circ p\right)^*\text{Hom}(E)$ at them. Getzler rescaling at a face $H$ consists of making use of the fact that the terms in the heat kernel that have order $k$ in the filtration (i.e. $k$-fold products of Clifford elements) come always with the factor $\rho_H^k$. Formally, one extends the filtrations off the faces via parallel transport on a transversal direction to them and shows that there exists a bundle\footnote{This bundle is a rescaling of the bundle $\left(\beta\circ p\right)^*\text{Hom}(E)$ in the same way $\leftindex^\phi TX$ is a rescaling of the bundle $TX$.} whose sections are precisely the sections in $\left(\beta\circ p\right)^*\text{Hom}(E)$ with such behavior close to the faces (à la Serre-Swan, see \cite[\S 8.1-8.4]{mel:aps}). Finally, one sets up a (rescaled heat) calculus (cf. \cite[p. 86]{vaillant}) with coefficients on the rescaled bundle and shows that the heat kernel belongs to it with the help of the Lichnerowicz formula for the square of the Dirac operator. This\footnote{ A similar philosophy has been used to prove an (less explicit) index statement in the (more general) setting of Lie manifolds \cite{bohlen-schrohe:getzler}.} is the invariant version of Getzler's coordinate computation \cite[p. 140-141, \S 4.3]{bgv}.

For example, if $c = 1$, we have the fibrations $\text{tf}\to X$ and $\phi\text{tf}\to B$ and the rescaled bundle is locally spanned by $\rho_{\text{tf}}\rho_{\phi\text{tf}}\frac{dx}{x}$, $ \rho_{\text{tf}}\rho_{\phi\text{tf}}dy_i$ and $\rho_{\text{tf}}xdz_j$ (see \cite[(87)]{vaillant}).

This is useful because the supertrace only sees terms with the highest Clifford degree (and kills off everything else), so we know that the leading order term recovered by the supertrace goes with the power of the boundary defining function of order the maximal Clifford degree in that face (in the example $c=1$, this allows us to compute the $n$-th order of the asymptotics at tf and the $(b+1)$-th at $\phi$tf).

We can rewrite the integral term for our $c$-$\phi$-metrics as:
$$\int_0^\infty \partial_t \leftindex^R {\text{Str }} e^{-t\dirac^2} dt = \int_0^\infty\int_{\partial X}  \left[\text{str}_p \left( \text{cl}\left(\frac{dx}{x^{2-c}}\right) \dirac e^{-t\dirac^2}\right)\right]_{b-c(b+f)}dy'dz'dt,$$

(with $f = \dim F$), meaning:
\begin{itemize}
    \item For $c = 0$, the contribution from $\dirac e^{-t\dirac^2}$ at the face $\phi$f we need is the coefficient of $\rho_{\phi\text{f}}^b$ ($\rho_{\phi\text{f}}$ being the boundary definition function for $\phi$f, think e.g.\footnote{The computations carried out within this method depend on the choice of boundary defining functions at the faces intersecting the lifted diagonal (because of the renormalized integrals in the McKean-Singer formula) and on the extension of the filtrations off the faces used to construct the rescaled heat bundle. However, the index formula obtained is independent of these choices.} $x'$). Since the heat kernel at this face is smooth and its expansion starts at order $0$ (\cite[Theorem 6.1]{tave} or Theorem \ref{thm:LongHeatKernelPoly}), its computation is not straightforward. Work in progress of the second author based on suggestions of P. Albin seems to indicate that a good approach is to use an idea of Bismut and Freed \cite[(2.6)]{bisfre} to relate the argument of the supertrace with the heat kernel of another Dirac operator squared and to perform Getzler rescaling in the $b$ base directions $dy_i$ and the normal direction $dx$. The rescaling \textit{localizes} and should thus yield a term of the form $\int_B \hat{\text{A}}(TB)\hat{\eta}(\dirac_F)$, where the Bismut-Cheeger $\eta$-form makes an appearance \cite[(0.5)]{bische}.
    \item For $c = 1$, Vaillant \cite[Theorem 4.11]{vaillant} showed that the asymptotics of $e^{-t\dirac^2}$ at bf start with $\rho_{\text{bf}}^{-f}$, which is precisely the contribution needed in the integral. As a result he obtains \cite[p. 97]{vaillant}:
    $$-\frac{1}{2}\eta(\dirac_{\mathcal{K},B}) = -\frac{1}{2\sqrt{\pi}}\int_0^\infty \text{Tr}_\mathcal{K}\left(\text{cl}\left(\frac{dx}{x}\right)\dirac_{\mathcal{K},B}e^{-t\dirac_{\mathcal{K},B}^2}\right)\frac{dt}{t^{1/2}}.$$
    \item For $c = 2$ there is no contribution (heuristically, there is no face in the heat blow-up space except at $t= 0$ or $t = \infty$: one should be able to show that the (renormalized) supertrace of the heat kernel is independent of time).
\end{itemize}

The face tf behaves in the same way for all metrics and one can use the classical arguments in it: the asymptotics start at $\rho_{\text{tf}}^{-n}$ (think $t^{-n/2}$) and thus we need to Getzler-rescale all $n$ directions to obtain the Atiyah-Singer integrand $\leftindex^R\int_X \hat{\text{A}}(\leftindex^{c - \phi}{TX})\text{Ch}(E)$ \footnote{Let $\nabla^E$ be the Clifford connection with which we construct our Dirac operator $\dirac = \text{cl} \circ \nabla^E$. Then, its curvature $\left(\nabla^E\right)^2$ is a differential form with values in $\text{hom}(E)\cong \mathbb{C}l(\leftindex^{c-\phi}TX)\otimes \text{hom}_{\mathbb{C}l}(E)$ and decomposes correspondingly in $R^E + \Omega^E$. The $\hat{\text{A}}$-genus is constructed à la Chern-Weil:
$$\hat{\text{A}}(\leftindex^{c-\phi}TX)=\text{det}^{1/2}\left(\frac{R^E/2}{\sinh \left(R^E/2\right)}\right).$$
Meanwhile, the Chern character is defined as a multiple of
$\text{Str}\left(e^{-\Omega^E}\right)$, so that no coefficients appear in front of its integral in the index formula.}.

For $c=1$, the face $\phi$tf is shown \cite[\S 5]{vaillant} to contribute with $-\int_B \hat{\text{A}}(TB)\hat{\eta}(\dirac_F)$ by a rescaling argument in the $dx$ and $dy_i$ directions. This is the same contribution one expects at $\phi$f (and perhaps at $\phi$tf for $c=2$).

Thus, the case $c = 1$ is completely understood by the thesis of Vaillant.

\begin{theorem}\cite{vaillant}
    Let $E\to X$ be a Clifford bundle over a manifold with a $1$-$\phi$-metric. Assume furthermore that the kernel of the vertical family $\dirac_F$ of the associated Dirac operator forms a vector bundle $\mathcal{K}\to B$. Then, the generalized index of $\dirac^+$ \cite[(67)]{vaillant}:
    $$\textup{ind}_-\left(\dirac^+\right) \coloneqq \lim_{t\to \infty}\leftindex^R{\textup{Str }}e^{-t\dirac^2} $$
    is given by\footnote{The different coefficients found in the literature for each of the terms correspond to different normalizations of the $\eta$, $\hat{\text{A}}$ and $\text{Ch}$ terms, see e.g. \cite[p. 21]{albroc} \cite[p. 47]{bgv}.}:
    $$\textup{ind}_- \left(\dirac^+\right) = \int_X \hat{\textup{A}}(\leftindex^{1-\phi}{TX})\textup{Ch}(E) - \int_B \hat{\textup{A}}(TB)\hat{\eta}(\dirac_F) -\frac{1}{2}\eta(\dirac_{\mathcal{K},B}).$$
    Moreover, $\dirac$ is Fredholm iff $\dirac_{\mathcal{K},B}$ is invertible and, in that case, $\textup{ind}_-\left(\dirac^+\right)$ corresponds to its Fredholm $L^2$-index.
\end{theorem}
\cite[\S 5.4]{vaillant} particularizes this result to the signature operator; the corresponding signature theorem is to be compared with \cite[Corollary 3]{hhm}, a work studying the spaces of $L^2$-harmonic forms for 0- and 1-$\phi$-metrics (generalized to QFB-metrics in \cite{KotRoc}).

The case $c = 0$ is taken up in the work of Leichtnam-Mazzeo-Piazza via a deformation and an adiabatic limit argument starting with Vaillant's result. This paper assumes full ellipticity and calls for a heat kernel approach to the problem, which we are currently developing.

\begin{theorem}\cite{lmp} \label{thm:index 0}
    Assuming the geometric Witt condition 
    $$\textup{Spec}(\dirac_F)\cap \left(-\delta, \delta\right) \quad \text{for some } \delta > 0,$$
    the Dirac operator associated to a Clifford bundle $E\to X$ over a manifold with a $0$-$\phi$-metric is fully elliptic and its index satisfies:
    $$\textup{ind}\left(\dirac^+\right) = \int_X \hat{\textup{A}}(\leftindex^{0-\phi}{TX})\textup{Ch}(E) - \int_B \hat{\textup{A}}(TB)\hat{\eta}(\dirac_F).$$
\end{theorem}

Finally, families index theorems for $c \in \{0, 1\}$ are also derived from Vaillant's
in \cite[Theorem 4.5, Corollary 5.1]{albroc}, where similar heat kernel rescalings take place.

We summarize the known results related to index theorems for $c-\phi$-metrics\footnote{Here, $\leftindex^{R}{\text{ind}}(\dirac^+)=\leftindex^{R}{\text{Str}}\left(\Pi_{\ker \dirac}\right)$. See also \cite[Theorem 4.25]{ash:thesis}.} for $c\in\{0,1,2\}$:

\scriptsize
\begin{table}\centering
\begin{tabular}{|l|p{35mm}|p{35mm}|p{45mm}|}
\hline
Metric & $0$-$\phi$ & $1$-$\phi$ & $2$-$\phi$ \\ \hline
\vspace*{\fill} Self-adjoint \vspace*{\fill} & \vspace*{\fill} \checkmark \vspace*{\fill} & \vspace*{\fill} \checkmark \vspace*{\fill} & Expected: unique closed extension if (\ref{eqn:geometric Witt}) \\ \hline
Fredholm & \vspace*{\fill} Iff $\dirac_F$ invertible \vspace*{\fill} & Iff $\dirac_{\mathcal{K},B}$ invertible, assuming $\mathcal{K} = \ker \dirac_F\to B$ vector bundle & \vspace*{\fill} Expected: always when self-adjoint, i.e. if (\ref{eqn:geometric Witt}) \vspace*{\fill} \\ \hline
\vspace*{\fill} $t\to 0$ \vspace*{\fill} & \vspace*{\fill} tf: $\int_X\hat{\text{A}}(\leftindex^{0-\phi}{TX})\text{Ch}(E)$ \vspace*{\fill} & \vspace*{\fill} tf: $\int_X\hat{\text{A}}(\leftindex^{1-\phi}{TX})\text{Ch}(E)$ \newline $\phi$tf: $- \int_B \hat{\text{A}}(TB)\hat{\eta}(\dirac_F)$ \vspace*{\fill} &  tf: $\leftindex^R\int_X\hat{\text{A}}(\leftindex^{2-\phi}{TX})\text{Ch}(E)$ \newline $\text{Expected } \phi\text{tf}$: $- \int_B \hat{\text{A}}(TB)\hat{\eta}(\dirac_F)$ \newline btf? \\ \hline
$t\in(0,\infty)$  & $\text{Expected } \phi\text{f}$: $- \int_B \hat{\text{A}}(TB)\hat{\eta}(\dirac_F)$  & bf: $-\frac{1}{2}\eta(\dirac_{\mathcal{K},B})$ & 0 \\ \hline
\vspace*{\fill} $t\to \infty$ \vspace*{\fill}  & $\phif_0$?
\newline zf: $\leftindex^{R}{\text{ind}}(\dirac^+)$  & \vspace*{\fill} $\text{ind}_-\left(\dirac^+\right)$ \vspace*{\fill} &  \vspace*{\fill} $\text{ind}\left(\dirac^+\right)$ \vspace*{\fill} \\ \hline
\end{tabular}
\end{table}

\normalsize
In the style of \cite{llarull}, one can use the index theory of $\dirac$ to look for obstructions to positive/non-negative scalar curvature metrics. This is based on the Lichnerowicz formula:
$$\dirac^2 = \left(\nabla^E\right)^*\nabla^E + \frac{\text{scal}_{g_c}}{4} + cl\left(\Omega^E\right).$$
As a result (and motivated by \cite{gromov}), there have been several advances in recent years on the topic of scalar curvature rigidity of spin manifolds. See \cite[\S 7]{albgr:incedge} for an application in a singular context similar to ours.

\section{Calderón projector and Dirichlet-Neumann operator}
\label{sec:calderon}

\begin{figure}
 \includegraphics[scale=1.1]{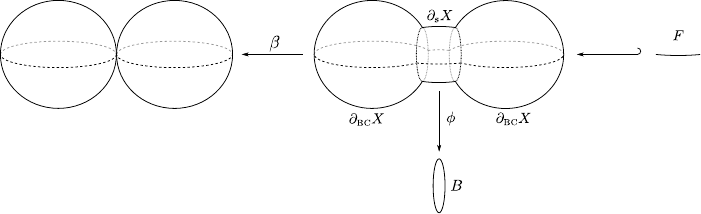}
 \caption{An example of a $\phi$-\textBC-manifold $X$ (exterior of the central picture) and how it arises via blow-up from a
 singular space (exterior of two touching spheres, left picture). Local coordinates are polar coordinates in the plane of tangency of the two spheres ($x=$ radius, $y=$ angle), and $z\in[-1,1]$ when parametrizing the left/right direction by $p(x)z$ where $p(x)\sim \const\, x^2$ as $x\to0$, with $z=\pm1$ corresponding to points on the spheres. Correspondingly, the resolution $X$ is obtained from the singular space in two steps: first, blow up the point of tangency of the two spheres; then the lifts of the spheres intersect the front face in a common circle $B=S^1$; next, blow up this circle. The second blow-up creates the
 singular boundary $\dsX$ as front face, so it is fibred over $B$. More precisely, $\dsX$ is the part of this front face which lies in the closure of the lift of the exterior domain, and the fibre is $F=[-1,1]$. The Euclidean metric on the exterior of the left picture lifts to a $2$-$\phi$-metric on $X$.
 }
 \label{fig:blow-up example}
\end{figure}

We now 
consider some aspects of boundary value problems for fibred cusp spaces with any $c\in\Z$. We need to extend our setting slightly. See Figure \ref{fig:blow-up example} for an illustration of the following definitions.

We define a $\phi$-\textBC-manifold, or $\phi$-\textbf{manifold with} \textBC-\textbf{boundary}, in a similar way as a $\phi$-manifold, except that the fibres may have non-empty boundary. More precisely, $X$ is now a manifold with corners
which has two types of boundary hypersurfaces: the `singular boundary hypersurfaces', given by $x=0$ in local coordinates and corresponding to $\bdX$ in Section \ref{sec:fib cusp}, and the `boundary hypersurfaces at which boundary conditions could be imposed', the `\textBC-boundary'. We assume that the boundary hypersurfaces of each type are pairwise disjoint, and denote the union of the singular ones by $\dsX$ and the union of the \textBC\ ones by $\dBCX$.
They intersect in their common boundary $\dsBCX$, which is a closed manifold. The  fibration \eqref{eqn:def phi-manifold} is replaced by a fibration
$\dsX \stackrel{\phi}\to B$, whose fibre $F$ now is a manifold with boundary, and 
the restriction of $\phi$ to $\partial(\dsX) = \dsBCX$, again denoted by $\phi$, defines a fibration $\partial F - \dsBCX \overset{\,\phi\,}\to B$.
Thus, the boundary $\partial F$ corresponds (near $x=0$) to $\dBCX$.
In local coordinates $x,y,z$, the \textBC-boundary is given by  $z_1=0$, and $z_1\geq0$ in $X$. For simplicity we assume again that $B$ is connected.

(General) $c-\phi$-metrics on $\phi$-\textBC-manifolds are defined as for $\phi$-manifolds, see \eqref{eqn:phi-forms}. In particular, they are non-degenerate up to the interior of the \textBC-boundary.
Here are some examples, which can also be combined (cf.\ Remark \ref{rem:connected}):
\begin{itemize}
 \item
 $c=2$: bounded domains in $\R^n$ with an incomplete cusp singularity (outward pointing) in the boundary, e.g.\ $\{(x,x^2z): x\geq0, |z|\leq 1\}$ near $x=0$. The Euclidean metric is roughly $dx^2 + x^4 g_F$ where $F$ is the $(n-1)$-dimensional unit ball, so $B$ is a point. $X$ is obtained by an iteration of two point blow-ups.
 
 Figure \ref{fig:blow-up example} shows an example with non-trivial base $B=S^1$.
 \item
$c=1$:  
domains in locally symmetric spaces (for example, take the central picture in Figure \ref{fig:examples}, but don't identify left and right side; instead, they form $\dBCX$)
 \item
 $c=0$: the slab $\R^{n-1}\times[0,1]$, or a uniformly fattened cone, for example an $\eps$-neighborhood of the right picture in Figure \ref{fig:examples} for small $\eps>0$. Another example is $\R^n\setminus\Omega$ for a bounded smooth domain $\Omega$. In either case we radially compactify Euclidean space at infinity (in the latter case $\dsX$ and $\dBCX$ are disjoint; this is really in the slightly extended setting of a disconnected base, see Remark \ref{rem:connected}; the results below extend to this setting, mutatis mutandis).
\end{itemize}

The notions of $\phi$-metric, {$\phi$-differential operator}, $\phi$-principal symbol, $\phi$-ellipticity, normal family and full ellipticity
carry over from the $\phi$-setting without \textBC-boundary verbatim, where we always assume smoothness up to $\dBCX$.
Also, $\dBCX$ with the restricted fibration is a $\phi$-manifold (without \textBC-boundary), so $\Psi^*_\phi(\dBCX)$ is defined.%

In \cite{FriGriSch:CPFCO} the authors consider Calderón projectors for the setting described above. 
Classically (i.e., in the non-singular setting $\dsX=\emptyset$), a Calder\'on projector associated with an  elliptic partial differential operator $P$ of order $m$ on a compact manifold  $X$ with non-empty boundary $\dBCX$ is a projection $C$ in $[\Cinf(\dBCX)]^m$ to the set $\calB_P$ of boundary data of solutions of the homogeneous equation
\begin{equation}
 \label{eqn:cauchy data space}
\calB_P = \{(u_{|\dBCX},\Dnu u_{|\dBCX},\dots, \Dnu^{m-1} u_{|\dBCX}):\, u\in\Cinf(X), Pu = 0\} \,.
\end{equation}
Here $\nu$ is some choice of transversal vector field in a neighborhood of the boundary and $\Dnu = \frac{1}{i}\partial_\nu$. Note there is no metric involved here.

It was first observed by Calder\'on \cite{Cal63} that such a projection exists which is a pseudodifferential operator, with an explicit principal symbol, and that this can be used to study boundary value problems. The first complete proof was given by Seeley in 1966, \cite{See66,See69}.
The result applies more generally to operators acting between sections of vector bundles.

When extending this result to our singular setting, i.e., to the case of a $\phi$-\textBC-manifold with $\dsX\neq\emptyset$, we need to specify 
the growth behavior at the singular set which we allow for the functions $u$ and for the vector field $\nu$ in \eqref{eqn:cauchy data space}. For $\nu$ we assume the geometrically natural condition that it is a $\phi$-vector field transversal to $\dBCX$.
For the functions $u$ there are various ways to restrict the growth behavior at $\dsX$. We use the letter $\calF$ to denote any choice of function space encoding such behavior which behaves well under restriction to $\dBCX$ and under the action of $\phi$-$\Psi$DOs on $\dBCX$. Thus, we have spaces $\calF(X)$, $\calF(\dBCX)$ of functions on $X$ and $\dBCX$, and similarly spaces of sections of vector bundles.
Typical examples are $\calF(X) = x^\alpha \Cinf(X)$, where $\alpha\in\R\cup\infty$, which would allow $x^{\alpha}$ decay (if $\alpha>0$) or growth (if $\alpha<0$) at the singular set, or spaces of functions which are conormal or polyhomogeneous at $\dsX$.
Then we define an  \textbf{$\calF$-Calder\'on projector} for $P$ to be a projection in $\calF(\dBCX)^m$ to the boundary data space \eqref{eqn:cauchy data space}, with $\Cinf(X)$ replaced by $\calF(X)$.
The main result of  \cite{FriGriSch:CPFCO}, stated with vector bundles,
is:
\begin{theorem} \cite{FriGriSch:CPFCO}
 \label{thm:Calderon}
 Let $X$ be a $\phi$-manifold with non-empty \textBC-boundary and $E,E'$ be complex vector bundles over $X$. Let $c\in\Z$ and
\begin{equation}
\label{eqn:cP operator}
P =x^{-cm}\tilde P,\ \tilde P\in\Diff_\phi^m(X;E,E') \,,
\end{equation}
 where $\tilde P$ is $\phi$-elliptic and satisfies the unique continuation assumption \ref{UCNF} stated below. Then there is an operator $C\in\Psi^*_\phi(\dBCX;E^m)$ which for any choice of admissible function space $\calF$  is an $\calF$-Calder\'on projector for $P$.
\\
When considered as acting between $m$-tuples of sections of $E$, the operator $C$ is an $m\times m$ matrix $(C_{kl})_{k,l=1\dots m}$,
where $C_{kl} \in \Psi_\phi^{k-l}(\dBCX;E)$.\footnote{\label{fn:connected Calderon}Note that in the example in Figure \ref{fig:blow-up example}, $\dBCX$, and in particular its boundary, $\dsBCX$, is disconnected. However, the base of the fibration $\dsBCX \overset{\,\phi\,}\to B$ is connected. By Remark \ref{rem:connected2} the $\phi$-double space of $\dBCX$ has $\phi$-front faces not only at, but also far away from the diagonal; this corresponds to non-trivial, i.e.\ not $O(x^\infty)$, behavior of 
the Schwartz kernel $K_C(p,p')$ of $C$ for $p$ on one of the spheres and $p'$ on the other, both approaching the singularity from the same direction.
This was stated slightly imprecisely in  \cite{FriGriSch:CPFCO}.
}

\end{theorem}

There is also precise information about the $\phi$-principal symbol and the normal family of $C$: both are Calderón projectors on suitable model spaces, for operators built from the $\phi$-principal symbol and the normal family of $P$. The full $\phi$-symbol of $C$ is determined constructively.

A Calderón projector is not unique since only its range is prescribed. There are ways to impose additional requirements to make it unique, for example by extending $P$ to an invertible operator on a manifold obtained by doubling $X$ across $\dBCX$, and taking the projection whose kernel is the boundary data space from the \lq other side\rq; or by choosing a metric on $\dBCX$ and on $E$ and using the projection which is orthogonal in the corresponding $L^2$ space.  This is discussed in detail in 
 \cite{FriGriSch:CPFCO}.

\begin{remark}
A particular issue in dealing with the Calder\'on projector is that $P$ may have \textbf{shadow solutions}, i.e.\ sections $u\not\equiv0$ satisfying $Pu=0$ whose boundary data at $\dBCX$ vanish. Ellipticity of $P$ then implies that they vanish to infinite order at $\dBCX$.
This is a type of failure of unique continuation for $P$. Similarly, an operator on the fibre $F$ may have shadow solutions (with respect to $\partial F$). Recall that the normal family is a family of operators on $F$. The assumption needed in Theorem \ref{thm:Calderon} is:
\begin{assumption}\label{UCNF} We assume that the normal families of
 $P$ and $P^\star$ have no shadow solutions.
\end{assumption}
It is not known if this assumption is necessary for the validity of the theorem.
It is shown in  \cite{FriGriSch:CPFCO} that any shadow solution of $Pu=0$ is rapidly decreasing at the singular boundary $\dsX$.
\end{remark}
The proof of Theorem \ref{thm:Calderon} proceeds along roughly the same lines as Seeley's proof in the non-singular case, involving extension of $P$ to an invertible operator $\Phat$ on a double of $X$, and  \cite{FriGriSch:CPFCO} provides a systematic exposition of this approach. The presence of the singular boundary poses additional challenges, in particular the invertibility of $\Phat$ necessitates invertibility of its normal family. This is the reason why Assumption \ref{UCNF} is needed. In addition, the \textbf{transmission property}, which is classical in the non-singular case, is  studied in detail for $\phi$-$\Psi$DOs, and rephrased in terms of conormal distributions, which may be of independent interest since it allows extension to other singular settings.

\medskip

In the forthcoming paper \cite{FriGriSch:DNOFCS} the results about the Calderón projector are used to study the Dirichlet-Neumann operator 
for a fibred cusp space. Recall that for a compact Riemannian manifold $(X,g)$ with non-empty boundary $\dBCX$ the \textbf{Dirichlet-Neumann operator} is the operator
\begin{equation}
\label{eqn:def DN op}
 \DN: \Cinf(\partial X)\to\Cinf(\partial X),\ f \mapsto \partial_\nu u\ \text{ where } \ \Delta u = 0,\ u_{|\partial X} = f \,,
\end{equation}
where $\partial_\nu$ is the outward unit normal derivative.
It is well-known that $\DN$ is a pseudo-differential operator of order 1, with principal symbol $|\xi|$. This can be proved using the fact that the Calderón projector is a pseudodifferential operator.

The main result of \cite{FriGriSch:DNOFCS} is that the Dirichlet-Neumann operator for a $c-\phi$-metric on a $\phi$-manifold with \textBC-boundary is well-defined and is in $x^{-c}\Psi^1_{\phi}(\dBCX)$. 
Note that if $g=x^{2c}g_\phi$ then 
$\Delta_{g} = x^{-2c} (\Delta_{g_\phi} + Q)$ where $Q$ is a first order $\phi$-operator. 
The extra factor $x^{-c}$ for $\calN$ arises from the normalization of the unit normal.
The result holds more generally if one replaces $\Delta_g$  by any elliptic operator $P\in x^{-2c}\Diff^2_\phi(X)$, under suitable injectivity assumptions, which are satisfied for $P=\Delta_g+\lambda$ if $\lambda\geq0$ and $c\geq0$, for example. 

In the case of a planar domain whose boundary has an incomplete cusp singularity $q$, the boundary is one-dimensional. Parametrizing one branch of it near $q$ by $x$, with $x=0$ at the singularity, the result means that the Schwartz kernel $K$ of $x^2\calN$ has the behavior explained in Figure \ref{fig:phi double space} as $x\to0$.\footnote{Also, if $x$, $x'$ correspond to points approaching $q$ on the two different branches, then the behavior of $K_{x^2\calN}(x,x')$ is the same but without singularity on the diagonal, i.e.\ as for operators in $\Psi^{-\infty}_\phi(\R_+)$. Compare Footnote \ref{fn:psido disconnected}.
}
In particular, $K(p,p')$ decays rapidly when $p,p'$ approach $q$ but have distance $d(p,p')>> x^2$. This remains true for general fibred incomplete cusps, for example in Figure \ref{fig:blow-up example} when  $p,p'$ approach the singularity from different directions, see Remark \ref{rem:X2phi geometry}.2.

The forthcoming paper \cite{GriSch:SDNOFCS} studies the spectrum of the Dirichlet-Neumann operator. The most notable result is that the spectrum is not discrete even in the incomplete case $c=2$, and a formula for the bottom of the essential spectrum. This was first observed in the special case of bounded domains with an incomplete cusp singularity (where the base is a point as noted above) in \cite{NazTas:SSPDWP}. For bounded domains with smooth boundary this spectrum (called the Steklov spectrum of the domain) is discrete since $\calN$ is an elliptic pseudodifferential operator of positive order on a closed manifold.

\section{Further directions, open problems}
\label{sec:more}

We list some problems involving fibred cusp spaces which are open  (as far as we know), as well as some research directions not treated here.

\begin{itemize}
\item
Instead of Dirac operators one may consider Dirac-Schrödinger operators, which in addition involve a potential. An index theorem for this case (\lq Callias index theorem\rq) was derived in the case of trivial fibre in \cite{kot:callias}, but the case of general fibre has not been treated.
\item Solve the question marks/items marked \lq expected\rq\ in the table after Theorem \ref{thm:index 0}.
\item
Apply the Calderón projector result to the study of boundary value problems (generalizing the classical theory, which can be found in \cite[Ch. 20]{Hor:ALPDOIII}).
\item
Study the Calderón projector and Dirichlet-Neumann operator for $c-\phi$-operators in the case where the base $B$ has boundary, not the fibre. A simple example for this is the Euclidean half space.
\item 
Study the local geometry of fibred cusp spaces near the boundary. The geodesics hitting a cuspidal singularity ($c=2$, base is a point) were analyzed in detail in \cite{GraGri:EMCS}: If $g$ is a $2$-$\phi$-metric, i.e.\ $g=x^4 g_\phi$ with $g_\phi$ as in \eqref{eqn:phi metric} without the $g_B$ term, then there is a unique geodesic starting at any point of $\bdX$, and they combine to a smooth exponential map based at the singular point. However, for typical (incomplete) embedded cuspidal singularities\footnote{A natural definition of these is:  a subset $M$ of a smooth Riemannian manifold has a cuspidal singularity at $q\in M$ if $M\setminus q$ is a submanifold and $M$ is resolved by blowing up $q$ and then a point on the front face. The smooth part of $M$ is
 equipped with the induced metric.} the metric does not have this form, but includes terms of the form $x^3 \,dx\,dz$,\footnote{More precisely, no matter how one chooses the coordinates, one cannot put it into this form.  This happens, for example, 
for surfaces similar to the left picture in Figure \ref{fig:examples}, but with non-circular link of the singularity. As shown in \cite{GraGri:EMCS} there is a geometric invariant, a certain function on $\partial X$, whose differential vanishes if and only if $g$ is a $2$-$\phi$-metric.} and this yields interesting unexpected behavior of the exponential map, for example it may be non-injective on any neighborhood of the singular point. The geodesics \textit{almost} hitting the singularity are studied in \cite{GriLye:GOS}, and in \cite{GriLye:WFGPTCN} this peculiar behavior is elucidated using a family of smooth spaces degenerating to the cusp singularity.

A similarly detailed study of the geodesic flow has not been carried out in the case where the base is non-trivial, or for most other singular geometries like (iterated) wedges.
\item
Develop the analogous spectral geometry in the case of iterated fibred cusp spaces, in particular studying their index theory and the appropriate notion of analytic torsion.
\end{itemize}

\section*{Appendix: Basics on manifolds with corners and blow-up}\label{sec:appendix}

We give a quick summary of basic notions on manifolds with corners. Details can be found in \cite{mel:aps}, \cite{Mel:DAMWC} or in the introductory text \cite{Gri:BBC}.

A \textbf{manifold with corners} of dimension $n$, denoted $X$ in the sequel, is defined like a manifold except that local charts are defined on open subsets of model spaces $\R^n_k:=\R^{n-k}\times\R_+^k$ for various $k\in\{0,\dots,n\}$ where $\R_+=[0,\infty)$, and an additional global condition is satisfied, see below. To define smoothness of transition maps or of maps between manifolds with corners we say that a map
 \[ \text{$U_1\to U_2$, where $U_i\subset\R^{n_i}_{k_i}$,} \]
is \textbf{smooth} if it extends to a smooth map
 \[ \text{$\Uwidetilde_1\to \R^{n_2}$, where $\Uwidetilde_1\subset\R^{n_1}$ is open and $U_1=\widetilde U_1\cap \R^{n_1}_{k_1}$.} \]
If $p\in X$ then there is a unique $k$, called the \textbf{codimension} of $p$, so that there is a coordinate system (inverse of a chart) mapping $p$ to $0\in\R^n_k$. The coordinates are then sometimes called \textbf{adapted} to $X$, and are often denoted $x=(x_1,\dots,x_k)$, $y=(y_1,\dots,y_{n-k})$ where $x_i\geq0$ and $y_j\in\R$ for all $i,j$.

A \textbf{face} of $X$ of codimension $k$ is the closure of a connected component of the set of points of codimension $k$. A \textbf{boundary hypersurface} is a face of codimension one. The global condition on a manifold with corners is that boundary hypersurfaces be embedded (rather than immersed) submanifolds. Equivalently, for each boundary hypersurface $H$ there is \textbf{boundary defining function} $\rho$ (often denoted $\rho_H$), i.e.\ a smooth function $\rho:X\to\R_+$ satisfying $\rho^{-1}(0)=H$ and $d\rho_{|p}\neq0$ for all $p\in H$.

A \textbf{p-submanifold} (where p is for product) of  $X$ is a subset $Y$ so that for each $p\in Y$ there is an adapted coordinate system on $X$ in which $Y$ is locally a coordinate subspace. That is, adapted coordinates $z=(x,y)$ can be chosen and  regrouped as $(z',z'')$ so that $Y=\{z''=0\}$ locally.
A connected p-submanifold $Y$ is called a \textbf{boundary p-submanifold} if $Y\subset\dX$, otherwise it is an \textbf{interior p-submanifold}. In the latter case only $y$ variables occur among the $z''$ variables.
For example, faces of $X$ are boundary p-submanifolds.

If $Y\subset X$ is a closed p-submanifold then the \textbf{blow-up} of $X$ in $Y$ is a new manifold with corners, denoted $[X,Y]$, together with a smooth map $\beta:[X;Y]\to X$, called \textbf{blow-down map}, which restricts to a diffeomorphism $[X;Y]\setminus\ff \to X\setminus Y$, where $\ff:=\beta^{-1}(Y)$ is called the \textbf{front face}, and so that
near any $p\in Y$ with coordinates $(z',z'')$ as above the map $\beta$ is locally near $\beta^{-1}(p)$ modelled by the polar coordinates map in the $z''$-coordinates, i.e.\ if $z'\in\R^{n'}_{k'}$, $z''\in\R^{n''}_{k''}$ then locally $\beta: \R^{n'}_{k'}\times \R_+\times \Sph^{n''-1}_{k''}\to\R^{n'}_{k'}\times\R^{n''}_{k''}, (z',r,\omega)\mapsto (z',r\omega)$ where $\Sph^{n''-1}_{k''}\subset\R^{n''}_{k''}$ is the unit sphere \lq octant\rq. Locally the front face is $\R^{n'}_{k'}\times\Sph^{n''-1}_{k''}$ and has local boundary defining function $r$. In practice, it is better to use \textbf{projective coordinates}. Blow-ups can also be iterated, i.e.\ if $Z$ is a p-submanifold of $[X,Y]$ then one can form $[[X;Y];Z]$ etc.

If $Z$ is a closed and connected subset of $X$ then the \textbf{lift} of $Z$ under the blow-up of a p-submanifold $Y\subset X$, denoted $\beta^*(Z)$, is defined as $\beta^{-1}(Z)$ if $Z\subset Y$ and as the closure of $\beta^{-1}(Z\setminus Y)$ if $Z\setminus Y$ is dense in $Z$, and undefined otherwise. If $Z$ is a p-submanifold meeting $Y$ \textbf{cleanly} (i.e.\ so that for each $p\in Y\cap Z$ there is an adapted coordinate system in which both $Y$ and $Z$ are coordinate subspaces) then $\beta^*(Z)$ is defined and a p-submanifold of $[X;Y]$. However, also subsets $Z$ which are not p-submanifolds can lift to such after (possibly iterated) blow-up, and then we say that $Z$ is \textbf{resolved} by the (iterated) blow-up. An important example is the diagonal $\{x=x'\}$ in $\R^2_+$, which is not a p-submanifold but is resolved by blowing up the origin. A p-submanifold $Y$ meets any face of $X$ cleanly, and the boundary hypersurfaces of $[X;Y]$ are $\ff$ and the lifts of the boundary hypersurfaces of $X$.

From now on we assume for simplicity that $X$ is compact.
The space of all smooth vector fields on $X$ which at each $p\in X$ are tangent to all boundary hypersurfaces containing $p$ is denoted $\calV_b(X)$.
Interpreting vector fields as first order differential operators and taking finite sums of smooth functions and compositions $V_1\circ\dots\circ V_l$ with all $V_i\in\calV_b(X)$, for $l\leq m$, we obtain the space of \textbf{b-differential operators} of order at most $m$, denoted $\Diff_b^m(X)$. Also $\Diff_b^*(X) := \bigcup_m \Diff_b^m(X)$. In adapted coordinates these are combinations of expressions $x_i\partial_{x_i}$, $\partial_{y_j}$ with smooth coefficients.

We now define various function spaces on a manifold with boundary $X$. The definitions can be extended to  manifolds with corners. First, $\Cinf(X)$ denotes  the space of functions $X\to\C$ which are smooth up to the boundary.
Let $x$ be a boundary defining function.
The following spaces consist of functions only defined and smooth on the interior of $X$ but having a certain prescribed behavior near the boundary:
\begin{itemize}
\item The space of \textbf{functions conormal to the boundary} of order $a\in\R$:
$$\calA^a(X)= \{u\in\Cinf(\interior{X}): \Diff_b^*(X)u \subset x^a L^\infty(X)\}\,,$$
and 
$\calA(X) = \bigcup_{a\in\R}\calA^a(X)\,.$
\item Spaces of functions \textbf{polyhomogeneous} at the boundary:
$$\calA^\calE(X) = \{u\in\calA(X):\, u\sim \sum_{(z,k)\in \calE} a_{z,k} \,x^z\log^k\!x\},
$$
where $\calE$ is an \textbf{index set}, i.e.\ a discrete subset of $\C\times\N_0$ satisfying certain additional conditions,  and $a_{z,k}\in\Cinf(X)$.
The symbol `$\sim$' means that the difference of $u$ and the finite sum over $\Re z\leq N$ is in $\calA^N(X)$, for each $N$. In particular, the asymptotics can be differentiated.
If $\calE=\N_0\times\{0\}$ then $\calA^\calE(X)=\Cinf(X)$. If the boundary is disconnected then each component can have its own index set (resp. order $a$ for $\calA^a(X)$). We write $\calE\geq a$ if $\calA^\calE(X)\subset\calA^a(X)$, i.e., $(z,k)\in\calE$ implies $\Re z\geq a$, and $k=0$ if $\Re z=a$.
\end{itemize}

For the definition of pseudodifferential operators we  need conormal distributions. A distribution $u$ on a manifold with corners $Z$ is \textbf{classical conormal} of order $m\in\R$ with respect to an interior p-submanifold $Y$ if it is smooth on $Z\setminus Y$ and near any point of $Y$ and in any adapted coordinate system $x,y=(y',y'')$ for $Z$ in terms of which $Y=\{y''=0\}$ locally,
\begin{equation}
 \label{eqn:def conormal}
u(x,y',y'') = \int e^{iy''\eta''}a(x,y';\eta'')\,d\eta''
\end{equation}
for a classical symbol $a$ of order $\mu=m+\frac14\dim Z - \frac12\codim Y$.
The space of these distributions is denoted $I^m_\cl(Z,Y)$.
We only need the case $\dim Z=2\dim Y$, then $\mu=m$.
Here classical means that $a$ has a symbol expansion
$a\sim \sum_{j=0}^\infty a_j$ where for each $j$
\begin{equation}
\label{eqn:pos homogeneous}
a_j (x,y';\lambda\eta'') = \lambda^{\mu-j} a_j(x,y';\eta'')
\end{equation}
for all $\lambda>0$ and all $x,y'$ and $\eta''\neq0$.
The $a_j$ are uniquely determined by $u$.

\bibliographystyle{amsalpha}
\bibliography{refs} 

\end{document}